\documentclass[11pt]{article}
\title{Extension of a theorem of Wschebor to free and matrix Brownian motions}

\usepackage{amsmath,amsthm,amscd,amssymb}
\usepackage{latexsym}
\usepackage{hyperref}
\usepackage{amsmath,amsthm,amscd,amssymb}
\usepackage{latexsym}
\usepackage[overload]{empheq}

\usepackage{bbm}
\usepackage{color}
\usepackage{ytableau}
\DeclareMathOperator{\GUE}{GUE}
\DeclareMathOperator{\GOE}{GOE}

\DeclareMathOperator{\SC}{SC}

\DeclareMathOperator{\cyc}{cyc}
\newcommand{\Cov}{\mbox{\rm Cov}}
\newcommand{\Var}{\mbox{\rm Var}}

\newcommand{\ve}{\varepsilon}
\newcommand{\tr}{\text{\rm{tr}}}

\newcommand{\W}{\mathcal W}

\newcommand{\R}{\mathbb R}
\newcommand{\He}{\mathbf H}
\newcommand{\I}{\mathrm{Id}_N}
\newcommand{\N}{\mathcal N}
\newcommand{\E}{\mathbb E}

\newcommand{\dd}{{\mathrm{d}}}

\newcommand{\vers}{\mathop{\longrightarrow}} 

\def \sur#1#2{\mathrel{\mathop{\kern 0pt#1}\limits^{#2}}}
\def \el{\sur{=}{(d)}}
\def\Tr{\hbox{Tr}\!\ }
\newcommand{\sn}{^{(N)}}
\newcounter{smalllist}

\numberwithin{equation}{section}
\newtheorem{theorem}{Theorem}[section]
\newtheorem{proposition}[theorem]{Proposition}
\newtheorem{lemma}[theorem]{Lemma}
\newtheorem{corollary}[theorem]{Corollary}
\theoremstyle{definition}
\newtheorem*{definition}{Definition}

\theoremstyle{remark}
\newtheorem{remark}{Remark}

\definecolor{Red}{rgb}{1,0,0}
\definecolor{Blue}{rgb}{0,0,1}

\usepackage{tikz}
\usetikzlibrary{decorations.pathreplacing}

\author{{\small Catherine Donati-Martin and Alain Rouault}
\footnote{Laboratoire de Math{\' e}matiques de Versailles, UVSQ, CNRS, Universit\'e Paris-Saclay, 78035-Versailles Cedex France. email : catherine.donati-martin@uvsq.fr, alain.rouault@uvsq.fr}
}
\begin{document}
\maketitle

\begin{abstract}
In 1992, M. Wschebor proved a theorem 
 on the  convergence of small increments of the Brownian motion.  Since then, it has been extended to various processes. We 
 prove a version of this theorem for  the Hermitian Brownian motion  and the free Brownian motion. Since these theorems deal with a convergence to a deterministic limit, 
 we prove also the convergence in distribution of the corresponding fluctuations. 
\end{abstract}
\bigskip

\section{Introduction}\label{intro}

In 1992 \cite{wschebor1992accroissements}, Mario Wschebor proved  the following
  remarkable property of the linear Brownian motion 
 $(W(t) , t\ge 0 \ ; W(0)=0)$.  If
\begin{align}
\label{W+}\mathcal W^\ve := \ve^{-1/2} 
\left(W(\cdot + \ve) -W(\cdot)\right)
\end{align}and if
 $\lambda$ is the Lebesgue measure on $[0,1]$, then,
 almost surely, for every $x \in \mathbb R$ and 
 every $t \in [0,1]$:
\begin{align}\label{W95}\lim_{\ve \to 0} \lambda \{s\le t:\, \mathcal W^\ve(s)  
\le x\} =  t \Phi(x)\,,   
\end{align}  
where $\Phi$ is the distribution function of the standard normal distribution $\mathcal N(0;1)$.
Motivated by the study of crossings of  stochastic processes in continuous time, he extended this result  to more general smooth approximations of trajectories. 

Let $(W_t, t \in \R)$ be the bilateral Brownian motion, i.e. $(W_t, t \geq 0)$ and $(W_{-t}, t \geq 0)$  are independent Brownian motions starting from $0$.
 \begin{theorem}[\cite{wschebor1992accroissements}, \cite{azais1996almost}, \cite{wschebor2006smoothing}]
 \label{W0}
  Let
  $\varphi : \R \to \R$ with compact support and bounded variation satisfying $\Vert\varphi\Vert_2 = 1$ and
 \begin{align}
 \label{def1.3}
W^\ve_\varphi (t)
=\int_{-\infty}^\infty \frac{1}{\sqrt \ve} \varphi\left(\frac{t-s}{\ve}\right) dW(s)\,,
 \end{align} 
and let
 \begin{align}
 \label{defmuve}
\mu^\ve = \int_0^1 \delta_{W^\ve_ \varphi(t)} dt\,,\end{align}
be the corresponding occupation measure.
Then almost surely as $\ve \to 0$, $\mu_\ve$ converges to the standard normal distribution $\mathcal N(0,1)$ 
 for the convergence of moments, i.e. that  for every $k \geq 0$
\begin{align}
\label{LLN}\int x^k d\mu^\ve (x) = \int_0^1 (W^\ve_\varphi(t))^k dt \vers_{\varepsilon \rightarrow 0} \E (\mathcal N^k)\,,
\end{align}
where $\mathcal N \el \mathcal N(0;1)$.
\end{theorem}
It is a statement which we call \textit{law of large numbers} (LLN). Actually (\ref{W+}) corresponds to $\varphi= 1_{[-1,0]}$ and the convergence of $\mu_\ve$ is equivalent to (\ref{W95}).

Since the Brownian motion $W$ is self-similar (Property P1) and has stationary increments (P2), it is possible to 
reduce the study of $\mu^\ve$ 
 ($\varepsilon \rightarrow 0$)
 to the study
 of the occupation measure in large time ($T :=\varepsilon^{-1} \rightarrow \infty$) for the rescaled process 
 \begin{align}
 \label{resc}W_\varphi^1(t) = \int_{-\infty}^\infty \varphi(t-s) dW(s)\end{align}
which does not depend  on $\varepsilon$.  This moving average  process is stationary Gaussian (it is the Slepian process when $\varphi= 1_{[-1,0]}$). Moreover, the independence of increments of $W$ (P3)  induces
 a finite dependence for $W_\varphi^1$ (recall that $\varphi$ is compactly supported), so $W_\varphi^1$ is ergodic, which allows to invoke  Birkhoff's theorem.
 
Later  the above result was extended to other types of processes (sharing properties P1 and P2), in particular the fractional Brownian motion (\cite{wschebor1995almost}, \cite{azais1996almost}, \cite{wschebor2006smoothing}). The Gaussian character of this process allows the use of its spectral measure instead of  independence of increments.

Corresponding fluctuations  (CLT) have been established  : for a large class of functions $g$ 
\begin{align}
\label{CLT1}
 \frac{1}{\sqrt\ve}
\left( \int_0^t g (W^\ve_\varphi(s)) ds - t \E g(\mathcal N)\right) \ , \ t\geq 0
\end{align}
converges, as $\varepsilon$ tends to zero to $(\sigma(g) W_t, t \geq 0)$ where $(W_t , t \geq 0)$ is a standard Brownian motion  and $\sigma(g)$ is an explicit function of $g$ (see \cite{berzin1997brownian} and \cite{wschebor2006smoothing} for a very complete review). 
 This follows from an 
 application of a continuous version of the  Breuer-Major theorem, which we recall now. 

Let 
$F$ be  a function from $\R$ to $\R$ such that $\E \left(F(\mathcal N)^2\right) < \infty$. It has the Hermite expansion
\begin{align}
\label{Hermiteexp}F(x) = \sum_{n=0}^\infty c_n H_n(x)\end{align} where $H_n$ is the Hermite polynomial of degree $n$. The Hermite rank of $F$ is  defined as 
the smallest $n$ such that the Hermite coefficient $c_n$ is nonzero. 

\begin{theorem}\cite{benhariz2002} \cite[Th. 1.1]{campese2020}
\label{theoscalfluc}
Let $X$ be a stationary Gaussian process with covariance function $\rho$. Assume that $F$ satisfies $\mathbb E F(\mathcal N)^p < \infty$ with $p\geq 2$  
and that its Hermite rank is $ \ell\geq 1$. Suppose also that $\int |\rho(t)|^\ell dt <\infty$.
Then the family of processes
\begin{align}
\label{2.3}
Z_\ve(t) := \ve^{1/2}\int_0^{t/\ve} F(X_s) ds \ , t\geq 0
\end{align} 
converges, as $ \ve$ tends to zero to $(\sigma_W(F) W_t, t \geq 0)$, where
\begin{itemize}
\item $(W_t, t \geq 0)$ is a standard Brownian motion,
\item $\sigma_W^2(F) = 2 \sum_{q= \ell}^\infty c_q^2 q! \int_0^\infty \rho(t)^q dt$\,, 
\item the convergence in law is fidi if $p=2$ and in $\mathcal C(\mathbb R^+)$ if $p > 2$.
\end{itemize}
\end{theorem}

The decomposition (\ref{Hermiteexp}) corresponds to the decomposition in chaos of the variable $F(X_s)$.

In the present paper, we extend these results  replacing the Brownian motion $W$  in (\ref{def1.3})
 by a  free Brownian motion and a matricial Brownian motion of dimension $N$, respectively.  

For the LLN, the deterministic limit is expressed in terms of the semi-circle distribution in the free case and in terms of the $\GUE(N)$ distribution in the matrix case. 

For the fluctuations in the free case, we use an extension of the Breuer-Major theorem due to \cite{kemp2012wigner} . The Hermite polynomials are replaced by the Chebyshev polynomials, and the Wiener chaos by the Wigner chaos.

In the case of matrix-valued fluctuations, no analogous chaos theory is available. We therefore require 
a suitable basis of orthogonal polynomials in the algebra of matrices. Recently,   Anshelevich and Buzinski \cite{anshelevich2022} introduced the notion of Hermite trace polynomials which is particularly adapted to our situation. These polynomials are indexed by a set of permutations. For a large class of these permutations we prove the convergence of fluctuations of the corresponding functional, to the law of an isotropic Gaussian matrix  (when we restrict to a fixed time $t$).

In Sec. \ref{sec:scal} we recall the scalar results, presenting short proofs to prepare the way for proofs in the free and matricial extensions. In Sec. \ref{sec:free} we establish the LLN and fluctuations for the free case, since the treatment is rather easy.  The more involved matrix case is addressed in Sections \ref{sec:matrix1} and \ref{sec:matrix2}. One of the main results of this paper concerns the matrix-valued  fluctuations in Theorem \ref{maintheoab}. 
The asymptotic regime $N\to \infty$ (from the matrix case to the free case) is summarized in diagram (\ref{diagram}) for the LLN and discussed in detail  in Section \ref{sec:asymp} for the fluctuations.
Besides, we consider in Sec. \ref{mvariate}  a different model of fluctuations using the Hermite matrix-variate polynomials. Finally, Section \ref{sec:app} is devoted to explicit computations needed in previous sections.
\section{The scalar case}
\label{sec:scal}
\subsection{LLN . Proof of Theorem \ref{W0}.}
\label{LLNs}
There are two proofs of the LLN, the historical one and the ``ergodic'' one. 
For the sake of completeness  we recall the proof of (\ref{LLN}) given by  Wschebor. 

 Let us assume that the support of $\varphi$ is contained in $[-a,a]$. The first step consists in stating that all moments of $\mu_{\mathcal W^\ve}$ converge in $L^2$ to the corresponding moments of $\mathcal N$. 
The marginals of $W_\varphi^\ve$ are  standard normal, 
 hence, for every $k$
\begin{align*}
\mathbb E \left(\int_0^1 W_\varphi^\ve(t)^k dt \right) = \int_0^1 \mathbb E  \left(W^\ve_\varphi (t)^k\right) dt  = \mathbb E\!\ (\mathcal N^k)\,.
\end{align*}
Besides
\begin{align*}\notag
\hbox{Var}\left(\int_0^1 W_\varphi^\ve(t)^k dt \right)&= \int_{[0,1]^2} \hbox{Cov}\left( W^\ve_\varphi(t)^k , W^\ve_\varphi(s)^k\right) ds dt\\
&= \ve^2\int_{[0,{1/\ve}]^2} \hbox{Cov}\left(W_\varphi^1(t)^k , W_\varphi^1(s)^k\right) ds dt
\,,
\end{align*}
where in the second identity we used the scaling property
\begin{align}
\label{el1}(W_\varphi^\ve(\ve t) , 0 \leq t\leq 1)  \el (W_\varphi^1(t), 0 \leq t \leq 1/\ve)\,.\end{align}
Then, we split the domain of integration into $\{(s,t) \in [0,1/\ve]^2 : |t-s| > 2a\}$ and ${[0,{1/\ve}]^2 \cap |t-s| \leq 2a}$. 
The first integral is zero since the process $W_\varphi^1$ is $a$-dependent, and the second integral is bounded, owing to Cauchy-Schwarz and
\[\lambda\{(s,t) \in [0,1/\ve]^2 : |t-s| \leq 2a\} \hbox{Var} (\mathcal N^k) =O(1/\ve)\,.\]
Then Borel-Cantelli and  properties of  Brownian paths allow to get an a.s. convergence. We  will not give details here.
\subsection{Fluctuations}
\label{sec:flucscal}
The fluctuations in the classical scalar case are a consequence of the continuous version of the Breuer-Major theorem \ref{theoscalfluc}.

\begin{corollary}
\label{cor1}
Let $\varphi$ satisfying the conditions of Theorem \ref{W0}. 
If $F$ is such that $\mathbb E F(\mathcal N^
2) < \infty$ and its Hermite rank is $\ell\geq 1$, then
\begin{align}
\label{cvflucscal}
\left(\ve^{-1/2}\int_0^t F(W_\varphi^\ve(s)) ds \ , t \geq 0\right) \vers_{\ve \to 0}^{(d)} \left(\sigma_W^2(F) W_t  \ , t \geq 0\right)
\end{align}
for the  fidi convergence, 
with 
 \begin{align}
 \label{defcov}
\rho(u) = \int_{-\infty}^\infty \varphi(u+\tau)\varphi(\tau) d\tau\,.
\end{align}
\end{corollary}
Notice that $\rho$ is defined on $\mathbb R$ compactly supported,  even and satisfy
\begin{align}
\label{CSrho}
|\rho(u)| \leq 1\,,
\end{align}
(use $\Vert \varphi\Vert_2 = 1$ and Cauchy-Schwarz). 
 All along the sequel we will use the quantity
\begin{align}
\label{defsigmadef}
\sigma_q^2 =  2 \int_0^\infty \rho(t)^q dt\,,
\end{align}
so that in the above (\ref{cvflucscal}) 
\begin{align}
\label{sigmaW}\sigma^2_W(F) = \sum_{q=\ell}^\infty q! c_q^2 \sigma_q^2\,.\end{align}

\begin{remark}
The above quantities $\sigma^2_q$ are non-negative, since they arise as limits of variances. If $F= H_1$ and $\int \varphi(s) ds = 0$, then $\sigma^2_F = \sigma^2_1 = 0$. In that case $\ve^{-1}\int_0^1 W_\varphi^\ve (s) ds $ has a Gaussian distribution whose variance does not depend on $\ve$.
\medskip

\noindent Examples :  

1) In Wschebor's original setting, 
\[\varphi=1_{]-1,0]}\ , \rho(t) = (1 - |t|)_+ \ ,   \int_0^\infty \rho(t)^q dt = 1/(q+1)\,.\]

2) If $\varphi = \frac{1}{\sqrt 2}\left(1_{[-1, 0[} - 1_{[0,1]}\right)$ we have
\[W_\varphi^\ve (s) = \frac{1}{\ve\sqrt 2}\left(W_{t+\ve} + W_{t-\ve} - 2 W_t\right) \ , \ \sigma_1^2 = 0\,,\]
and $\ve^{-1}\int_0^1 W_\varphi^\ve (s) ds \el \mathcal N(0; 2/3)$.
\end{remark}
\medskip

 In  \cite{BLL-book} there are two proofs of the convergence of one-dimensional marginals in Th. \ref{theoscalfluc} :  a classical one and a modern one based on the Fourth Moment Theorem. Both are based on a spectral representation of the process.  
 
 We give here a variation of the proof of   Corollary \ref{cor1},
  starting from the moving average representation itself and multiple Wiener-Itô integrals.
  
  We will begin with a fixed chaos, i.e. $F = H_n$, a fixed time $t=1$ in  (\ref{2.3}) and $X_s = W_\varphi^1(s)$. We set $\ve = T^{-1}$ . 
We have
\begin{align}
\label{defintH}
T^{-1/2}\int_0^T H_n\left(W_\varphi^1(t)\right) dt = I_n^W(f_T)\end{align}
where
\begin{align}
\label{deffT} f_T(t_1, \dots, t_n) = T^{-1/2}\int_0^T \varphi(t-t_1) \cdots \varphi(t-t_n) dt\,,\end{align}
and for  $f$ symmetric function in  $L^2(\R_+^n)$
\begin{align}
\label{defIW}
I_n^W(f) = 
\int_{\R_+^n} f(s_1, \dots, s_n) dW_{s_1} \cdots dW_{s_n}\,.
\end{align}
This yields
\begin{align*}
\notag
\E \left(I_n^W(f_T)\right)^2 &= n! \Vert f_T\Vert^2 = n! T^{-1}\int \left(\int_0^T \varphi(t-t_1) \cdots \varphi(t-t_n) dt\right)^2 dt_1 \cdots dt_n\\
\notag
&= n! T^{-1}\int_{[0,T]^2}\left(\int \varphi(t-t_1)\varphi(s-t_1) dt_1\right)^n dt ds\\
&=  n! T^{-1} \int_{[0,T]^2} \rho(t-s)^n dt ds \,,
\end{align*}
where we used Fubini.
After a change of variable, we conclude that
\begin{align}
\notag
\E \left(I_n^W(f_T)\right)^2 = 2n!  \int_0^T \left(1- \frac{u}{T}\right) \rho(u)^n du\\
\label{EIpB2}
 \vers_{T\to \infty} 2 n! \int_0^\infty \rho(u)^n du = n! \sigma_n^2\,.
 \end{align}
 We apply  the Fourth Moment Theorem (\cite[Th. 1]{Nualart2005Pec}) which 
 says that the convergence in distribution of $I_n^W(f_T)$ to a normal variable as $ T \to \infty$ is equivalent to the convergence of the fourth moments and also equivalent to
  the convergence in $L^2$ to zero of contractions
 \begin{equation}
 \notag
f_T\otimes_k f_T (\xi_1, \dots, \xi_{2n-2k}) =
\end{equation}
\begin{equation*}
= \int f_T(s_1, \dots, s_k, \xi_1,\dots, \xi_{n-k})f_T(s_1, \dots s_k, \xi_{n-k+1},\dots, \xi_{2n-2k}) ds_1 \dots ds_k \,,
 \end{equation*}
 for $k=1,\dots , n-1$. 
 But 
 \begin{align}
 \notag
 \Vert  f_T\otimes_k f_T \Vert^2 = T^{-2}\int_{[0,T]^4}\rho(t-\tau)^k \rho(s-\sigma)^k\rho(t-s)^{n-k}\rho(\tau -\sigma)^{n-k} dt ds d\tau d\sigma\,,
 \end{align}
 which tends to $0$ (see p.11-12 in \cite{BLL-book}). This proves 
 \begin{align}
 \label{BM}
 I_n^W(f_T) \vers_{T \to \infty}^{(d)} \mathcal N(0, \sigma^2_n)
 \,.
\end{align}
\noindent In view of proving the fidi convergence of the process
\[T^{-1/2}\int_0^{Tu} H_n\left(W_\varphi^1(t)\right) dt \ , \ u \geq 0\]
to the Brownian motion, we apply \cite[Th. 1]{peccati2005gaussian}. We have just to check that if for $a<b$ 
\[f_T^{[a,b]}(s_1, \dots, s_n) := \frac{1}{\sqrt T}\int_{aT}^{bT}\varphi(t-s_1) \cdots\varphi(t-s_n) dt\]
then for $u_1 < u_2 \leq u_3 < u_4$
\begin{align}
\label{covW}\E\left(I_n^W(f_T^{[u_1, u_2]})I_n^W(f_T^{[u_3, u_4]})\right) \vers_{T \to \infty} 0\end{align}
but this expectation  is equal to
\begin{align*}
\notag
n! \langle f_T^{[u_1, u_2]}, f_T^{[u_3, u_4]}\rangle_{L^2} = \frac{n!}{T}\int_{[Tu_1, Tu_2]\times[Tu_3, Tu_4]}\rho(t-\tau)^n dt d\tau\\
\notag
=  \frac{n!}{T}\int_{\tau \in [Tu_1, Tu_2], v+\tau \in [Tu_3, Tu_4]} \rho(v)^n dvd\tau\\
=  n! \int \lambda\left([u_1,u_2]\cap [u_3- vT^{-1}, u_4 - vT^{-1}]\right) \rho(v)^n dv\,,
\end{align*}
where $\lambda$ is the Lebesgue measure. It is straightforward to see that the limit is $0$.
 This proves the fidi convergence of increments.
 
This argument can be directly extended  to obtain  the convergence for   combinations of Hermite polynomials. 
Eventually, tightness is proved in \cite{campese2020}.
\begin{remark}
The above results can also be extended to moving average processes driven by more general processes sharing with Brownian motion properties P1 and P2 (see Sec. 1), such as Hermite processes (see  
\cite{nualart2021}).
\end{remark}
\section{The free case}
\label{sec:free}
A non-commutative probability space is an algebra $\mathcal A$ of operators  on a complex separable Hilbert space, closed under adjoint and convergence in the weak operator topology, equipped with a trace $\tau$, 
that is, a weakly continuous linear functional, satisfying:
\begin{itemize}
\item $\tau(\mathbf 1) = 1$
\item $\tau(ab) = \tau(ba)$
\item $\tau(aa^*) \geq 0$ and $\tau(aa^*) = 0$ iff $a=0$\,.
\end{itemize}
The self-adjoint elements of $\mathcal A$ are called non-commutative random variables. If $a$ is such a random variable, 
 the linear form on the set of  polynomials
  defined by
\[P \in \mathbb C[X] \mapsto \tau(P(a))\,,\]
is called the distribution of $a$. In this case there exists a 
 unique probability measure $\mu$ such that 
\[\tau(P(a)) = \int_\mathbb R P(x) d\mu(x) \,.\]
The semi-circle distribution 
 of variance $\sigma^2$ is 
\[\SC (0; \sigma^2)(dx) =  \frac{1}{2\pi\sigma^2}\sqrt{(4\sigma^2 -x^2)_+}\,\ dx\,.\]
If $a$ is a non-commutative $\SC(0;1)$ ($\SC$ for short)  random variable, its moments are given by
\begin{align}
\label{Cat}\tau(a^{2p}) = \frac{1}{p+1}\binom{2p}{p}= C_p \ , \tau(a^{2p+1}) = 0\,,\end{align}
where $C_p$ is the $p$th Catalan number. In particular $\tau(a^2) = 1$ and $\tau(a^4) =2$.

The non commutative variables $a_1, a_2, \dots, a_p$ are called free if one has for all $n$, $i_1, i_2, \ldots i_n \in \{1, \ldots, p\}$,
$$ \tau (a_{i_1} a_{i_2} \ldots a_{i_n}) = 0$$
when $\tau(a_j) = 0$ for $j = 1,  \ldots,  p$ and $i_1 \not= i_2 \not= \ldots \not=i_n$ (consecutive indices are distinct).

The free Brownian motion is a family of  
$(\mathbf S(t), t \geq 0)$ of elements of $\mathcal A$ such that:
\begin{itemize}
\item $\tau(\mathbf S (t)) = 0$ for all $t$
\item $\mathbf S$ has free increments : for $t_1 < t_2 < \ldots t_p$, $\left((S(t_{i+1}) - S(t_i) , i=1, \ldots, p-1\right)$ are free.
\item $\mathbf S(0) = 0$ and the distribution of $\mathbf S(t) - \mathbf S(s)$ is $\SC(0;t-s)$ for $s< t$.
\end{itemize}
As in Sec. \ref{intro}, we can also define  a bilateral free Brownian motion, i.e. $(\mathbf S(t) , t \geq 0)$ and $(\mathbf S(-t), t \geq 0)$ are free Brownian motions  which are mutually free.

Such a  process is self-similar with index $1/2$. Now we define
\begin{align}
\label{defSbfve}
S_\varphi^\varepsilon (t) =  \ve^{-1/2} \int_{-\infty}^\infty \varphi\left(\frac{t-s}{\ve}\right) d\mathbf S(s)\,,
\end{align}
where as usual, we assume that the kernel $\varphi$ is bounded, has a support included in $[-a, a]$  and satisfies $\Vert \varphi\Vert_2 =1$.  In all that follows, $P$ will be a polynomial in one variable.
\subsection{LLN}
\begin{proposition}
\label{LLNfree} For the convergence of moments, we have:
\begin{align}
\label{limjk}
\lim_\ve \int_0^1 P(\mathbf S_\varphi^\ve(t)) dt = \int_0^1 P(x) d\SC(x) \cdot \mathbf 1 \,.
\end{align}
\end{proposition}
\proof 
It is clear that we may suppose that $P(x) =x^k$.  We have to prove that for  every $j$
\begin{align}
\label{4b}
\lim_{\ve\rightarrow 0}\tau\left(\int_0^1 \left(\mathbf S_\varphi^\ve (t)\right)^k dt\right)^j = (m_k)^j\,,
\end{align}
where  $m_k= \tau(a^k)$ when $a$ is a $\SC$ random variable.

 We adapt the proof of the scalar case Sec. \ref{LLNs}. 
The basic properties are scaling  :
\begin{align}
\label{elfree1}
(\mathbf S_\varphi^\ve(\ve t) , 0 \leq t\leq 1)  &\el (\mathbf S_\varphi^1(t), 0 \leq t \leq 1/\ve)\,,
\end{align}
and stationarity, in particular 
\begin{align}
 \label{elfree2}
 \mathbf S_\varphi^\ve (t) &\el \SC\,.
 \end{align}
Let us first prove (\ref{4b}) for $j=1$ and $j=2$. For $j=1$ there is no limit to take, since from (\ref{elfree2})  
\[\tau \left(\int_0^1 \left(\mathbf S_\varphi^\ve(t)\right)^k \ dt \right) = \int_0^1 \tau\left(\left(\mathbf S_\varphi^\ve(s)\right)^k\right)ds = m_k\,.\]
For $j=2$, it is enough to prove 	
\begin{align}
\label{j=2}\lim_\ve \tau\left(\int_0^1\left(\mathbf S_\varphi^\ve(t)^k - \tau\left(\mathbf S_\varphi^\ve(t)^k\right)\right) dt\right)^2 = 0\,.\end{align}
We have successively
 \begin{align}\notag
& \tau\left(\int_0^1\left(\mathbf S_\varphi^\ve(t)^k  - 
 \tau\left(\mathbf S_\varphi^\ve(t)^k\right)\right) dt\right)^2 =\\
 \notag
 &= \tau\left(\int_{[0,1]^2} 
 \left(\mathbf S_\varphi^\ve(t)^k - \tau\left(\mathbf S_\varphi^\ve(t)^k\right)\right)
 \left(\mathbf S_\varphi^\ve(s)^k - \tau\left(\mathbf S_\varphi^\ve(s)^k\right)\right)dt ds\right)
\\
\notag
&= \int_{[0,1]^2}  \tau\left(\left(\mathbf S_\varphi^\ve(t)^k - \tau\left(\mathbf S_\varphi^\ve(t)^k\right)\right)
 \left(\mathbf S_\varphi^\ve(s)^k - \tau\left(\mathbf S_\varphi^\ve(s)^k\right)\right)\right)dt ds\\
 \notag
 &=\ve^2 \int_{[0,1/\ve]^2}  \tau\left(\left(\mathbf S_\varphi^1(t)^k - \tau\left(\mathbf S_\varphi^1(t)^k\right)\right)
 \left(\mathbf S_\varphi^1(s)^k - \tau\left(\mathbf S_\varphi^1(s)^k\right)\right)\right)dt ds\,. 
 \end{align}

 As in the scalar case we split $[0,1/\ve]^2$ 
  into $\{(s,t) \in [0,1/\ve]^2 : |t-s| > 4a\}$ and ${[0,{1/\ve}]^2 \cap |t-s| \leq 4a}$. 
The first integral is zero since  $\mathbf S_\varphi^1(t)^k$ and $\mathbf S_\varphi^1(s)^k$ are free as soon as  $|t-s| > 4a$,

 In the second integral we have by Cauchy-Schwarz
\begin{align}
\notag| \tau\left(\left(\mathbf S_\varphi^1(t)^k - \tau\left(\mathbf S_\varphi^1(t)^k\right)\right)
 \left(\mathbf S_\varphi^1(s)^k - \tau\left(\mathbf S_\varphi^1(s)^k\right)\right)\right)|\leq
 \end{align}
 \begin{align}
 \notag \leq \sqrt{\tau\left(\mathbf S_\varphi^1(t)^k - \tau\left(\mathbf S_\varphi^1(t)^k\right)^2\right)}\sqrt{\tau\left(\mathbf S_\varphi^1(s)^k - \tau\left(\mathbf S_\varphi^1(s)^k\right)^2\right)}
 = m_{2k}- (m_k)^2\,,
 \end{align}
 so that the second integral is
\[O\left(\lambda\{(s, t)\in [0, 1/\ve]^2  : |t-s| \leq 4a\}\right) =O(1/\ve) \,,\]
and then, after multiplication by $\ve^2$ its contribution vanishes.

Assume $j\geq 3$. Setting
\[\mathcal J_j^\ve = \left(\int_0^1 \left(\mathbf S_\varphi^\ve (s)\right)^k ds\right)^j\,, \]
we have
\begin{align*}
\mathcal J_j^\ve - (m_k)^j = (\mathcal J_1^\ve - m_k) \Delta_j^\ve\,,
\end{align*}
where $
\Delta_j^\ve = \Sigma_{r=0}^{j-1} \mathcal J_r^\ve (m_k)^{j-r}$ 
and by Cauchy-Schwarz again
\begin{align}
\label{CS}
\left|\tau (\mathcal J_j^\ve) - (m_k)^j\right|\leq \sqrt{\tau(\mathcal J_1^\ve -m_k)^2}\sqrt{ \tau((\Delta_j^\ve)^2)} 
\end{align}
Now since the stationary law $\SC$ has a compact support $[-2,2]$, we have for every $p$
\begin{align*}
\Vert \mathcal J_p^\ve \Vert \leq 2^{kp}\,,
\end{align*}
so that, for every $\ve$ 
\begin{align*}
\tau((\Delta_j^\ve)^2) \leq j^2  2^{2kj}
\end{align*}
and then
\[\tau (\mathcal J_j^\ve) - (m_k)^j = O\left(\sqrt{\tau\left(\mathcal J_1^\ve  -m_k\right)^2}\right)
\]
which tends to zero by (\ref{j=2}). $\Box$
\subsection{Fluctuations}
\label{flucfree}
Let $U_n$ be the $n$th Chebyshev polynomial of the second kind, defined on $[-2,2]$ by
\begin{align}
\label{defChe}
U_n (2 \cos \theta) = \frac{\sin ((n+ 1)\theta)}{\sin \theta}\,.
\end{align}
These are the monic orthogonal polynomials associated to the $\SC$ distribution.
\begin{proposition}
\label{prop:freefluc}
Let $P  = \sum_{q=1}^d c_q U_q$  be a polynomial  of degree $d$ centered for the $\SC$ distribution, then 
\begin{align*}
\left(\ve^{-1/2}\int_0^t P(\mathbf S_\varphi^\ve (s)) 
 ds\right)_{t\geq 0}  
\ \vers_{\varepsilon \rightarrow 0}^{(d)}
\  \left(\sigma_S(P)\mathbf S(t)\right)_{t\geq 0} 
\end{align*}
where
\begin{align*}
\sigma_S(P)^2 = \sum_{q=1}^d  c_q^2 \sigma_q^2\,,
\end{align*}
where  $\sigma_q$ is defined in  (\ref{defsigmadef}) 
and the convergence is fidi.
\end{proposition}
We recall some facts about the convergence  of  free Wigner chaos associated with the  free Brownian motion  $(\mathbf S_t)_{t\geq 0}$ (we refer to Biane and Speicher \cite{biane1998stochastic} for stochastic calculus for free Brownian motion).
 
For $f \in  L^2 (\R_+^n)$ with some symmetry, one can define a multiple stochastic integral of order $n$ with respect to $\mathbf S$, denoted by $I_n^S (f)$ (see  \cite{kemp2012wigner}). In particular, for $f = 1_{[0,t]}^{\otimes n}$, 
$I_n^S(f) = U_n(\mathbf S_t)$. \\

A striking result in \cite{kemp2012wigner} gives a criterion for the convergence of noncommutative variables in   a fixed Wigner chaos. 
Let $(f_k)$ be a sequence of $L^2$ functions in $\R^n_+$. In \cite{kemp2012wigner}, the authors proved that the convergence in distribution of a sequence $(I_n^S(f_k))_k$ as $k \rightarrow \infty$ to a semicircular distribution is equivalent to the convergence of the  fourth moment to 2 (for a normalized sequence satisfying $||f_k||_2 = 1$) and is also equivalent to a condition expressed in terms of contractions of the $(f_k)$. \\
As a corollary, they obtain the following Wiener-Wigner transfer principle.

\begin{theorem}\cite[Th. 1.8]{kemp2012wigner}
\label{transfer}
Assume that the functions $f_k$ are fully symmetric. Let $\sigma$ be a finite constant. Then, as $k\to \infty$,
\begin{enumerate}
\item
$ \E[I_n^W(f_k)^2] \to n!\sigma^2$ if and only if $ \E[I_n^S(f_k)^2]\to \sigma^2$\,, 
\item
If the relations above  are verified, then 
\[I_n^W(f_k) \vers^{(d)} \mathcal N(0 ; n! \sigma^2)\] if and only if 
\[I_n^S(f_k) \vers^{(d)} \SC(0;\sigma^2)\,.\]
\end{enumerate}
\end{theorem}

\noindent\textbf{Proof of Prop. \ref{flucfree}}
Let us first consider a fixed index $n$ and a fixed time $t=1$. We 
 want to prove 
 the following result (analogous to the Breuer-Major theorem) : 
\begin{equation} \label{SC}
\frac{1}{\sqrt{\varepsilon}} \int_0^1 U_n( \mathbf S_\varphi^\varepsilon ( s) )ds \  \vers_{\varepsilon \rightarrow 0}^{(d)} \SC(0; \sigma_n^2) 
 \,.
\end{equation}
 By scaling property,
$$\frac{1}{\sqrt{\varepsilon}} \int_0^1 U_n( \mathbf S_\varphi^\ve ( s) )ds  \stackrel{(d)}{=}  \sqrt{\ve} \int_0^{\frac{1}{\ve} } U_n(\mathbf S_\varphi^1(s)) ds$$
and,  with $\ve^{-1}= T$, like (\ref{defIW})
\begin{align}
\label{defintHS}
T^{-1/2}\int_0^T U_n\left(\mathbf S_\varphi^1(t)\right) dt = I_n^S(f_T)
\end{align}
with $f_T$ defined in (\ref{deffT}) and for $f \in L^2(\R_+^n)$ 
\begin{align}
I_n^S(f) = \int_{\R_+^n} f(s_1, \dots, s_n) d\mathbf S_{s_1} \cdots d\mathbf S_{s_n}\,.
\end{align}

We can now apply the Wiener-Wigner transfer principle Theorem \ref{transfer} 
  to assert that from the convergence \eqref{BM}, we deduce the convergence in distribution \eqref{SC}. 
 
  Let us extend the above result to a polynomial $P(x) = \sum_1^d c_k U_k(x)$. We consider
  \[Y_T =T^{-1/2}\int_0^T P(\mathbf S_\varphi^\ve(s)) ds = \sum_1^d c_k Y_{k, T}\]
  where $Y_ {k, T}$ is associated with $U_k$. We know that $Y_{k, T}$ converges in distribution to $\SC(0 ; \sigma_k^2)$. Moreover,  from the multi-dimensional Wiener-Wigner transfer theorem \cite[Th. 1.6]{nourdin2013multi}, the vector $(Y_{1, T}, \dots, Y_{d, T})$ converges in distribution to $(S_1, \dots, S_d)$ where the $S_i$ are $\SC(0 ; \sigma_i^2)$ distributed and free, since they belong to  distinct chaoses. To show that $Y_T$ converges to $\SC(0; \sum_1^d c_k^2\sigma^2_k)$, we have to consider the convergence of moments. Actually
  \[(Y_T)^j = \left(\sum_1^d c_k Y_{k, T}\right)^j= Q(Y_{1, T}, \dots, Y_{d, T})\] 
for some polynomial $Q$. From the above convergence, we deduce the convergence of $ Q(Y_{1, T}, \dots, Y_{d, T})$ to $Q(S_1, \dots,  S_d)=(\sum_1^d c_k S_k)^j$. This means that $Y_T$ converges in  distribution to $\SC(0; \sum_1^d c_k^2 \sigma_k^2)$.

Let us fix $n$.  In view of proving the fidi convergence of the process
\[T^{-1/2}\int_0^{Tu}U_n\left(\mathbf S_\varphi^1(t)\right) dt \ , \ u \geq 0\]
to the free Brownian motion, we refer to Sec. \ref{sec:flucscal}. 
We consider four times $u_1 <u_2 <u_3< u_4$ and the limit (\ref{covW}). Applying again 
the multi-dimensional Wiener-Wigner transfer theorem \cite[Th. 1.6]{nourdin2013multi} we deduce that 
\[\left(I_n^S(f_T^{[u_1,u_2]}) , I_n^S(f_T^{[u_3,u_4]})\right) \vers_{T \to \infty}^{(d)} (S_1, S_2)\]
where $S_1, S_2$ are free and distributed as $\SC(0; (u_2-u_1), \SC(0,; u_3-u_4))$.

 An extension to a general polynomial and the fidi convergence can be obtained by combining the above arguments. It is left to the reader.
  $\Box$
\medskip

In particular for $\varphi= 1_{]-1, 0]}$ we obtain
\begin{corollary}
\label{cor:Un}
\begin{align*}
\frac{1}{\sqrt \varepsilon}\int_0^1 U_n\left(\frac{S(s+\varepsilon) - S(s)}{\sqrt \varepsilon}\right) ds  \vers_{\varepsilon \rightarrow 0}^{(d)} \SC\left(0 ; \frac{2}{n+1}\right)\,.
\end{align*}
\end{corollary}
\section{Matricial case : LLN}
\label{sec:matrix1}
In this section, we consider the space $\mathcal H_N$ of $N\times N$ Hermitian matrices, and we denote by $\tr$ the normalized trace $N^{-1}\Tr$. 
We replace the bilateral Brownian motion of Sec. \ref{sec:scal} with a bilateral  Hermitian Brownian motion 
 $W\sn$,  with values in the space $\mathcal H_N$ 
  such that 
$(W\sn_{i,j}, i\leq j)$ are independent bilateral Brownian motions,  complex if $i < j$  and real if $i=j$ with variance \[\mathbb E |W\sn_{i,j}(t)|^2 =|t|/N\,.\]
Set
\[\mathbb W_\varphi^{\ve, N}(t) =\ve^{-1/2}\int_{-\infty}^\infty \varphi\left(\frac{t-s}{\ve}\right) d W\sn(s)\,,\]
with the usual assumption: $\varphi$ compactly supported and satisfying $\Vert\varphi\Vert_2 = 1$.
For $t=1$, the distribution is  the classical Gaussian Unitary Ensemble of size $N$ and variance $N^{-1}$, which we will denote by $\GUE(N^{-1})$ in the sequel. We denote by $M\sn$ a random matrix distributed as $\GUE(N^{-1})$.
We will often omit the index $\varphi$ when it will be clear.
\begin{proposition}
\label{matversp} Assume that $N$ is fixed. Then, 
almost surely
\begin{align}
\label{matvers}
\forall k  \geq 1 \ ,\ \lim_\ve \int_0^1 \left(\mathbb W^{\ve, N}(t)\right)^k dt = \mathbb E (M\sn)^k\,.\end{align}
\end{proposition}
\proof 
 The proof is componentwise. We have
 \begin{align}
 \notag
 \left[\int_0^1 \left(\mathbb W^{\ve, N}(t)\right)^k dt\right]_{i,j} =  \sum_{i_2, \dots, i_{k-1}} 
 \int_0^1 \left(\mathbb W^{\ve, N}(t)_{i, i_2}\dots \mathbb W^{\ve, N}(t)_{i_{k-1},j}\right) dt\,,
 \end{align}
Each integrand is selfsimilar, $2a$-dependent  and stationary as in Sec. 1, so 
\begin{align*}
\lim_\varepsilon \int_0^1 \left(\mathbb W^{\ve, N}(t)_{i, i_2}\dots 
\mathbb W^{\ve, N}(t)_{i_{k-1},j}\right) dt = \mathbb E\left(M\sn_{i, i_2}\dots M\sn_{i_{k-1},j}\right)
\end{align*}
a.s. and we get (\ref{matvers}) by summation.  $\Box$
\bigskip

\noindent\underline{Let us now consider the asymptotics $N\to\infty$.}
\medskip

In all the sequel, $P$ will be a polynomial in one variable. Let us consider the following diagram

\begin{equation}
\label{diagram}
\begin{CD}
\int_0^1  P\left(\mathbb W^{\ve, N}(s)\right)  ds \in \mathcal H_N @>(1): \ve \rightarrow 0>>  \mathbb E \left( P(M\sn)\right) \in \mathcal H_N\\ 
@V(3) : N \rightarrow \infty VV @V(2) : N\rightarrow \infty VV\\
\int_0^1 P\left(\mathbf S^\ve (s)\right) ds\in \mathcal A  @>(4) : \ve \rightarrow 0>> \int_0^1 P(x) d\SC(x) \cdot \mathbf 1 \in \mathcal A\,.
\end{CD}
\end{equation}
\medskip

Let us give the precise meaning of all these arrows, at least when $P$ is a monomial $P(x) = x^k$. 
(1) is the result of Prop. \ref{matversp}   and (4) was proved in Sec. \ref{LLNfree}.
  
  Let us prove (2). Set $H_N := \mathbb E \left( P(M\sn)\right)$ and $c :=  \int_0^1 P(x) d\SC(x)$.
(2)
  means that 
 \begin{align*}
  \lim_{N\rightarrow \infty} \tr H_N^j = c^j, \ \hbox{for all} \ j\,.
  \end{align*}
  Actually, 
  \[H_N^j = \E\left[P(M\sn_{(1)})\cdots P(M\sn_{(j)})\right]\]
  where $M\sn_{(1)},\dots, M\sn_{(j)}$ are i.i.d. $\GUE(N^{-1})$.  
From the Wigner's theorem, we have then
\[\lim_{N\rightarrow \infty} \tr H_N^j  = \tau \left[P(S_1) \cdots P(S_j)\right]\,,\]
where $S_1, \dots, S_j$ are free and $\SC$.  The RHS is clearly $c^j$.

Let us prove (3). 
Assume  $\ve = 1$ to simplify notations.
Set
\[J^N_j =  \tr \left(\int_0^1  \mathbb W^{1,N}(s)^k ds\right)^j\,.\]

We have to prove that
\begin{align*}\lim_N \mathbb E J_j^N = \tau\left(\int_0^1 (\mathbf S^1(s)^k ds\right)^j
\end{align*}
and, if possible
\begin{align}
\label{possible'}\lim_N  J^j_N =\tau\left(\int_0^1 (\mathbf S^1(s)^k ds\right)^j \ \hbox{(a.s.)}\,.\end{align}
\underline{Let us  begin with $j=1$.}

Since $\mathbb W^{1, N}(s) \el M^{(N)}$ and 
$\mathbf S^1(s)  \el \SC$,
we have
\begin{align*}
\notag
\E J^N_1 &= 
 \tr \int_0^1 \E\left(\mathbb W^{1,N}(s)^k\right) =  \tr\ \E  (M^{(N)})^k\\
\tau \left(\int_0^1 \mathbf S^1(s)^k\right) ds &= \int_0^1 \tau\left(\mathbf S^1(s)^k\right)ds = m_k\,,
\end{align*}
hence from the convergence (2)
\begin{align*}
\lim_N \E J^N_1 = \lim_N \E \tr \left(M^{(N)}\right)^k  = m_k =
\tau \left(\int_0^1 \left(\mathcal S^1(s)\right)^k\right) ds \,.
\end{align*}
To prove the a.s. convergence for $j=1$, we will use the second moment method.
 Let us know look for $ \Var J^N_1$. From the classical Fubini-like formula
\[\Var \left(\int_0^1 T(s) ds\right) = \int_{[0,1]^2} \Cov(T(s), T(u)) ds du\,,\]
where $T$ is any $L^2$ process, 
we get
\[ \Var J^N_1 = \int_{[0,1]^2} \Cov\left(\tr \mathbb W^{1,N}(s)^k, \tr \mathbb W^{1,N}(u)^k\right) ds du\,.\]
Now, by Cauchy inequality and stationarity
\begin{align*}
\left|\Cov\left(\tr \mathbb W^{1,N}(s)^k,  \tr \mathbb W^{1,N}(u)^k\right) \right| \leq \Var \left(\tr (M^{(N)})^k\right)\,, 
\end{align*}
which tends to $0$ as $N \rightarrow 0$ (see \cite[Sec. 2.1.4]{agz}.
Actually, since the bound is $O(N^{-2})$, the convergence may be strengthened to an a.s. convergence, owing to Borel-Cantelli lemma.

\underline{For $j \geq 2$}, 
let us rewrite $J_j^N$ as a multiple integral
\begin{align*}
J_j^N = \tr \int_{[0,1]^j} \mathbb W^{1,N}(s_1)^k \cdots \mathbb W^{1,N}(s_j)^k ds_1 \cdots ds_j\,,
\end{align*}
so that, by Fubini
\begin{align}
\label{Fubini}
\mathbb E J_j^N = \int_{[0,1]^j} \mathbb E  \tr \left( \mathbb W^{1,N}(s_1)^k \cdots \mathbb W^{1,N}(s_j)^k \right) ds_1 \cdots ds_j\,,
\end{align}
It should be clear that the process $\mathbb W^{1,N}$ converges to $\mathbf S^1$ in the sense that for every sequence $(s_1, \dots, s_p)$
\[\lim_N \  \mathbb E \tr [\mathbb W^{1,N}(t_1) \cdots\mathbb W^{1,N}(t_p)] = \tau \left(\mathbf S^1(t_1) \cdots\mathbf S^1(t_p)\right)\,.\]

So, we have a convergence pointwise of the integrand in (\ref{Fubini}) to 
\[\tau\left(\mathbf S^1(s_1)^k \cdots \mathbf S^1(s_j)^k\right)\,.\]
Up to  an application of the Lebesgue dominated theorem, we could conclude (again by Fubini)
\begin{align*}
\lim_N \mathbb E J_j^N &=  \int_{[0,1]^j}\tau\left(\mathbf S^1(s_1)^k \cdots \mathbf S^1(s_j)^k\right) ds_1 \cdots ds_j\\
&= \tau \left(\int_{[0,1]^j}\mathcal S^1(s_1)^k \cdots \mathbf S^1(s_j)^k ds_1 \cdots ds_j\right)
= \tau \left(\int_{[0,1]} \mathbf S^1(s)^k ds\right)^j\,.
\end{align*}
Now, we have to find a uniform bound for $\mathbb E  \tr \left( \mathbb W^{1,N}(s_1)^k \cdots \mathbb W^{1,N}(s_j)^k \right)$. 
By Hölder inequality
\[|\mathbb E \tr (A_1\dots A_j)| \leq  \prod_{r=1}^j \left(\mathbb E  \tr (A_r)^j\right)^{1/j}\,,\]
so that
\begin{align}
\notag
|\mathbb E \tr \left( \mathbb W^{1,N}(s_1)^k \cdots \mathbb W^{1,N}(s_j)^k \right)|\leq \prod_{r=1}^j \left(\mathbb E \tr(\mathbb W^{1,N}(s_r)^{kj}\right)^{1/j}
&=  \mathbb E \tr(M_N)^{kj}
\end{align}
which is bounded in $N$  since it converges to $\tau \left((\mathbf S^1)^{kj}\right)$.

\section{Matricial case : fluctuations}
\label{sec:matrix2}

As seen in the previous sections,  Hermite polynomials in the scalar case and  Chebyshev polynomials in the free case play a major role in the study of fluctuations.  In the matrix case, we will use Hermite trace polynomials, which are matrix-valued functions with matrix argument. They are defined in the next subsection.  This can be regarded as an intermediate link between the scalar and free cases. 

 Section \ref{sec:fluc} is devoted to the main result and its proof. Finally, in Section \ref{sec:asymp} we let $N \to \infty$ and recover some results in terms of free probability.


\subsection{Hermite polynomials for Hermitian Brownian motion}
This section summarizes the ``stochastic part''  of the paper by Anshelevich and Buzinski \cite{anshelevich2022} in which the authors define  Hermite trace polynomials for the Hermitian Brownian motion. \\
First, let us recall what is known for the real Brownian motion $(B_t)_{t\geq 0}$. There exists polynomials $H_n(x,t)$ such that for all $n$, $H_n(B_t, t)$ is a martingale for the filtration induced by $B$. More precisely, $H_n(x,t) = t^{n/2} H_n(\frac{x}{\sqrt{t}})$ where $H_n$ is the classical Hermite polynomial of degree $n$. Moreover, we have the chaos representation
$$ H_n(B_t,t) = \int_{[0,t]^n } dB_{t_1} \ldots dB_{t_n}.$$

Let $(W^{(N)}(t))$ be a $N \times N$ Hermitian Brownian motion as in Sec. 4. We cannot find a polynomial  $P(x,t)$ of degree $n$ in $x$ such that $P(W^{(N)}(t),t)$ is a martingale, for $n \geq 3$. For example, for $n =3$, a martingale $M^{(N)}$ involving $(W^{(N)}(t))^3$ is given by :
\begin{equation} \label{Herm3}
M^{(N)} (t) = (W^{(N)}(t))^3 - t( 2W^{(N)}(t) +\tr(W^{(N)}(t)))\,. 
\end{equation}
Note that for $N=1$, we recover the classical Hermite polynomial $H_3(x,t) = x^3 -3t x$. 

Therefore, we need to replace the class of polynomial with a larger class of trace polynomials. Informally, a trace polynomial in a matrix indeterminate $X$  is a linear combination of product of the form $X^k \prod_i \tr(X^{l_i})$.

From now in this section, we drop  the superscript $N$ for simplicity and write $W_t$ instead of $W\sn(t)$.
 For $h \in L^2(\R)$, we set 
\[W(h) = \int_{-\infty}^\infty h(t) dW_t\,,\]
and clearly $W_t =W(h)$ for $h=1_{[0,t]}$.\\
In \cite{anshelevich2022}, the authors define trace polynomials as follows. Let $S_0(n)$ denote the set of permutations of $[0,n]:= \{0,1, \ldots, n \}$. Given  a collection $(A_i)$ of $N\times N$ matrices and  $\alpha \in S_0(n)$  define
 \begin{align}
\notag
\Tr_\alpha (A_1, \dots, A_n) &=\\
\label{parttrace} \prod_{ \textrm{ $i$  in the cycle starting with 0} } &A_i\prod_{ \textrm{ others  cycles }} \Tr \left (\prod_{j \textrm{ in the cycle }} A_j\right).
\end{align}
The matricial $\alpha$-Hermite polynomial $H_\alpha (W(h))$ is defined in (\ref{defHermite}) as linear combination of trace polynomial $\mathrm{Tr}_\eta (W(h), \ldots, W(h))$, where $\eta$ is a contraction $C_\pi (\alpha)$ of $\alpha$, see 
 \cite[Def. 3.13,  Th. 6.4, Rem. 6.5]{anshelevich2022}. We give some examples of Hermite polynomials in Prop. \ref{4H}. We recall the  definition of these contractions below.\\
\begin{definition}
\label{defiAns} (see \cite[Def. 3.9]{anshelevich2022})
Let $\pi$ be a partition of $[n]:= \{1, \ldots, n\}$ into blocks of size 1 or 2, i.e. $\pi \in \mathcal{P}_{1,2}(n)$. Denote by $\textrm{supp}(\pi) = [n] \backslash \textrm{Sing}(\pi)$ the set of indices belonging to  2-blocks. For $\alpha \in S_0(n)$ and $\pi \in \mathcal{P}_{1,2}(n)$,  define
\begin{align}\notag
 C_\pi(\alpha) &= \left(\frac{1}{N}\right)^{a(\alpha, \pi)} \beta_\pi(\alpha)\\
 \label{def:Cpi}
\beta_\pi(\alpha) &= P^{[0,n] \backslash \textrm{supp}(\pi)} _{[0,n-2l]} (\pi \alpha) \big|_{\textrm{supp}(\pi)^c} \,,
\end{align}
where \begin{enumerate}
\item $l$ is the number of 2-blocks of $\pi$.
\item For $A, B$ two sets of integers with the same cardinal, $P^A_B$ is the unique order-preserving bijection from $A$ to $B$ and the corresponding bijection on the set of permutations on $A$, resp. on $B$. 
\footnote{Example (\cite[Ex. 3.8]{anshelevich2022}): For $\alpha =(13524)$ and $S= \{2,5\}$, $\alpha|_{S^c} = (134)$ and $P_{[3)}^{[5]\setminus S}\alpha|_{S^c} =(123)$.}

\item $a(\alpha, \pi) = \textrm{cyc}_0((\pi \alpha)|_{\textrm{supp}(\pi)^c} -\textrm{cyc}_0(\pi \alpha) + l$ where $\textrm{cyc}_0$ denotes the number of cycles of a permutation on $ [0, n]$, not containing 0\,. 
\end{enumerate}
Finally, for $\alpha\in S_0(n), u \in \R$ and $M \in \mathcal H_N$, set
\begin{align}\label{2var}
\widetilde H_\alpha(M, u) = \sum_{\pi \in \mathcal{P}_{1,2}(n)} (-1)^{n-|\pi|}  \left(\frac{1}{N}\right)^{a(\alpha, \pi)} u^l\mathrm{Tr}_{\beta_\pi(\alpha)}\left(M, \dots, M\right)
\end{align}
where the trace polynomial $\Tr$ has $n-2l$ arguments.

\noindent The $\alpha$-Hermite polynomial is defined by 
\begin{equation} \label{defHermite}
H_\alpha (W(h))= \widetilde H_\alpha (W(h), \Vert h\Vert^2)\,,
\end{equation}
\end{definition}
\begin{remark} 
 Note that if $\pi$ has only 2-cycles, $C_\pi(\alpha)$ is equal to $(0)$ up to multiplicative constant. 
\end{remark}
We present some explicit examples in the following proposition
\begin{proposition}
\label{4H}Let  $\alpha_n =(01\dots n)$.
The first $\widetilde H_{\alpha_n}$ are :
\begin{align*}
\widetilde H_{\alpha_1} (M, u) &= - M\,,\\
\widetilde H_{\alpha_2}(M, u &= M^2 -u\,,\\
\widetilde H_{\alpha_3}(M, u)&= -M^3 + 2u M +u \tr M\,,\\
\widetilde H_{\alpha_4}(M,u)&= M^4 -2 u M\tr M- 3u M^2 -u \tr M^2 + (2+ N^{-2})u^2\,.
\end{align*}
\end{proposition}

When $N=1$,  definition \eqref{defHermite} recovers the classical Hermite polynomials (up to a sign).  
This discrepancy in sign arises from the choice of normalization adopted in \cite{anshelevich2022}, which 
is based on  algebraic motivations  rather than probabilistic conventions. The proof of this proposition is given in Section \ref{sec:app}.
\medskip

In the next  two propositions, we present some  moment formulas for tracial polynomials of the Hermitian Brownian and for products of Hermite polynomials. Their notations are required for our main Theorem \ref{maintheoab} and their statements are used in its proof.
   
\begin{proposition} See  \cite[Prop. 6.2]{anshelevich2022} Let $\{D_1, \ldots D_n\}$ be 
non-random Hermitian matrices. Then for even $n$, for $\alpha \in S_0(n)$, 
\begin{align}
\notag
\E(\mathrm{Tr}_\alpha  (W(h_1)D_1, \ldots, W(h_n)D_n) ) =&\\
 \label{mixmom} \frac{1}{N^{n/2}} \sum_{\pi \in \mathcal{P}_2(n)} C_\pi( h_1 \otimes \ldots \otimes h_n) &\mathrm{Tr}_{\pi \alpha} (D_1, \ldots, D_n).
\end{align}
where $\mathcal{P}_2(n)$ is the set of pair partitions of $[n]$ and
$$C_\pi(h_1 \otimes \ldots \otimes h_n) = \prod_{(i,j) \textrm{ pair of } \pi} \langle h_i, h_j\rangle.$$
\end{proposition}
The following proposition comes from \cite[Prop. 5.16]{anshelevich2022} and the identification at the beginning of Sec. 5 therein.
\begin{proposition}
\label{prop:moment} 
 For $\alpha_n = (01\ldots n)$, $k \in \mathbb{N}$ and $\alpha_{nk} = (01\ldots nk)$,
\begin{eqnarray}\notag
\E( H_{\alpha_n}( W(h_1)) \ldots H_{\alpha_n}( W(h_k)) ) =&\\
\notag
\sum_{\pi \in \mathcal{P}_2(n, \ldots, n)} C_\pi(\alpha_{nk}) \ &C_\pi( h_1^{\otimes n} \otimes \ldots \otimes h_k^{\otimes n}) =\\
\label{moment}
 \sum_{\pi \in \mathcal{P}_2(n, \ldots, n)}  \left(\frac{1}{N}\right)^{ nk/2 - \textrm{cyc}_0(\pi \alpha)} &C_\pi( h_1^{\otimes n} \otimes \ldots \otimes h_k^{\otimes n})\,.
\end{eqnarray}
where $\mathcal{P}_2(n, \ldots, n)$ is the set of inhomogeneous pair partitions of $[nk]$, meaning that if $(i,j)$ is a pairing of the partition with $i<j$, then $i \leq pn < j$ for some $1 \leq p \leq k-1$ (i.e. considering the set $nk$ as $k$ blocks of $n$ elements, any pairing involves two elements of different blocks). \\
In \eqref{moment}, with a slight abuse of notation, $C_\pi(\alpha_{nk})$ denotes the constant  K in the writing of the partition $C_\pi(\alpha_{nk}) = K (0)$ and the LHS is equal to the RHS multiplied by the identity matrix. 
\end{proposition}
\subsection{Fluctuations }
\label{sec:fluc}
Let   $(W_t, \ t\in (-\infty, \infty))$ be  as above. 
Let $\varphi : \R \to \R$ compactly supported with bounded variation and satisfying   $\Vert\varphi\Vert_2 = 1$. Set $\varphi_t(s) = \varphi(t-s)$.
We now consider the stationary process with value in  $\mathcal{H}_N$  defined by : 
$$X_t = \int_{-\infty}^\infty \varphi (t-s) dW_s := \int_{-\infty}^\infty \varphi_t(s) dW_s = W(\varphi_t)\,.$$
Notice that 
\begin{align}
\label{hh}\langle \varphi_t, \varphi_s\rangle = \rho(t-s)\end{align} where
$\rho (t) = \int \varphi(t+u) \varphi(u) du$   as in \ref{defcov}. \\
Let $H_{\alpha_n}$ the Hermite polynomial on  $\mathcal{H}_N$ with $\alpha_n = (01\ldots n) \in S_0(n)$. \\
As in the scalar case, we are looking for a CLT for the random matrix :
$$M_T = \frac{1}{\sqrt{T}} \int_0^T  H_{\alpha_n} (X_s) ds$$ as $T$ tends to $\infty.$
The main result of this section is :
\begin{theorem}
\label{maintheoab}
As $T \rightarrow \infty$, $M_T$ converges in distribution to a Gaussian random matrix $M_\infty$ whose law is characterized by the following properties :
\begin{enumerate}
\item[1)] It is invariant by conjugation which implies that  in the decompostion $M_\infty = U^* D U$, the matrix $U$ of eigenvectors is Haar distributed and independent of the matrix $D$ of eigenvalues.
\item[2)] The law of the spectral distribution 
\[\mu_{M_\infty}  = \frac{1}{N} \sum_{i=1}^N \delta_{\lambda_i(M_\infty)}\] is given by its moments : for $k$ even,
\begin{equation}
\E\left(\int x^k d\mu_{M_\infty} (x)\right) = \E(\tr (M_\infty^k)) = K_{k, n,N} \sigma_n^k
\label{momentlim}
\end{equation}
 where $\sigma^2_n$ is defined in (\ref{defsigmadef}),
 \begin{itemize}
 \item
  $K_{k, n,N} = \sum_{\pi \in \Pi_{n,k}} C_\pi(\alpha_{nk})$  where 
\[C_\pi(\alpha_{nk}) =   \left(\frac{1}{N}\right)^{ nk/2 - \textrm{cyc}_0(\pi \alpha_{nk})}\]
(see  (\ref{def:Cpi}))
  \item $\Pi_{n,k}$ is the set of partitions $\pi$ of $\mathcal{P}_2(n, \ldots, n)$ in inhomogeneous pairs such that if  $i_1$ and $i_2$ are in the same block of size $n$, then $\pi(i_1)$ and $\pi(i_2)$ are also in the same block.
  \end{itemize}
\item[3)] We have 
\begin{align}
\label{second}
M_\infty \el \sigma_n\left(a_{n,N} G + b_{n,N} \xi I_N\right)
\end{align}
where
\begin{itemize}
\item
$G \el \GUE(N^{-1})$, $\xi\el \N(0;N^{-2})$ (real) and $G$ and $\xi$ are independent,
\item
the coefficients $a_{n,N}$ and $b_{n,N}$ are defined via the permutation $\tilde\alpha = (0)(1\dots n)(n+1\dots 2n)$ by 
 \begin{align}
 \label{defa}
 a^2_{n,N} &= \frac{1}{N^n}\sum_{\pi \in \mathcal P'_2(n,n)} N^{ \cyc (\pi\tilde\alpha)}\\
 \label{defb}
b^2_{n,N}&=\frac{1}{N^n}\sum_{\pi \in \mathcal P_2(n,n) \setminus \mathcal P'_2(n,n)} N^{ \cyc (\pi\tilde\alpha)}
 \end{align}
 with 
 \begin{align*}
 \mathcal P'_2(n,n) = \{\pi \in \mathcal P_2(n,n) : n \ \hbox{and} \ 2n \ \hbox{are in the same cycle of} \  \pi\tilde\alpha\}\,.
\end{align*}
\end{itemize}
 \end{enumerate}
\end{theorem}
For example, we provide the values of $a_{n,N}$ and $b_{n,N}$ for $n=2$ and $n=3$ (see Section \ref{sec:app} for details).

\begin{proposition}
\label{prop:ab3}
 For $n=2,3$, the coefficients in the decomposition (\ref{second}) are 
\begin{align}
\notag
a_{2,N} &= b_{2, N} = 1\,,\\
\label{ab3}
a_{3, N}^2 &= 1+3N^{-2} \ , \ b_{3,N}^2 =  2\,. 
\end{align}
\end{proposition}
\begin{remark}
We choose to study the permutation $\alpha_n$ since the leading term of $H_{\alpha_n}$ is $X^n$. It is obviously restrictive, but given the length of the computations, we preferred not to make the paper any  longer. Of course, studying of fluctuations for a general functional of the process would require
\begin{itemize}
\item a complete catalogue of limits when $\alpha_n$ is replaced by other permutations
\item the use of one of the chaos decomposition formula in Anshelevich.
\end{itemize}
We also choosed to restrict ourselves to a fixed time.
\end{remark}
\medskip

\begin{remark}
It may seem natural to tackle this problem with the tools of Wiener chaos, as in the scalar case. 
Unfortunately, we did not find any ready-made version of chaoses applicable to matrix-valued Brownian motion.
It is noteworthy that in  \cite[Remark 1.40]{kemp2012wigner}, the authors leave this ``involved'' task  to the reader.
\end{remark}

 \proof  The convergence of $M_T$ to a Gaussian random matrix follows from Breuer-Major CLT theorem for non-linear functional of stationary Gaussian process (extended to multidimensional Gaussian process in \cite{arcones}). Indeed, for any non-random Hermitian matrix $A$, we can write
$$\Tr(AM_T) =\frac{1}{\sqrt{T}} \int_0^T F(X_s) ds$$ for some function $F$ from $\mathcal{H}_N$ to $\R$, centered for the multidimensional stationary Gaussian process $(X_t)$, and thus converges to a  centered  Gaussian random variable of variance $\sigma^2 (A, n, N)$. It follows that $M_T$ converges to a Gaussian Hermitian matrix.

1)  Since the distribution of $(X_t)_{t \geq 0}$ is invariant by unitary conjugation, it follows that for all unitary matrix $U$, $M_T \stackrel{(d)}{=} UM_TU^*$ and therefore $M_\infty \stackrel{(d)}{=} U M_\infty U^*$.  
We deduce (see \cite[Cor. 2.5.4]{agz}) that the eigenvectors of $M_\infty$ are Haar distributed and are independent of the eigenvalues. 

2) Let us  now compute, for $k \geq 1$, the moment
\[\int x^k d\mu_{M_\infty} (x) =\E(\tr (M_\infty^k))\,,\]
and prove (\ref{momentlim}). We have 
\begin{align*}
M_T^k =
T^{-k/2}\int_{[0, T]^k} H_\alpha(X_{t_1})H_\alpha(X_{t_2})\dots H_\alpha(X_{t_k}) dt_1\cdots dt_k\,.
\end{align*}
From \eqref{moment},
this yields:
\begin{align*}
\mathbb E M_T^k =
 \sum_{\pi \in \mathcal{P}_2(n, \ldots, n)} C_\pi(\alpha_{nk}) \  T^{-k/2}\int_{[0, T]^k}  C_\pi( \varphi_{t_1}^{\otimes n} \otimes \ldots \otimes \varphi_{t_k}^{\otimes n}) dt_1\cdots dt_k\,.
\end{align*}
From (\ref{hh})  
 we have
\begin{align*}
C_\pi \left(\varphi_{t_1}^{\otimes n} \otimes \cdots\otimes \varphi_{t_k}^{\otimes n}\right) 
= \prod_{1 \leq i < j \leq k} \rho(t_i - t_j)^{n_{ij}}
 \end{align*}
where   $n_{ij}$ denotes the number of pairs consisting of an element of  the block $i$ and an element of  block $j$.  Of course  $0 \leq n_{ij}\leq n$. If $n_{ij} = n$ the blocks $i$ and $j$ are  said to be completely connected. 
Let us denote by $J$ those pairs of blocks such that $n_{ij} =n$ and let  $|J|$ be its cardinal. We can split 
 \begin{align}
 \notag
 \int_{[0,T]^k} C_\pi \left(h_{t_1}^{\otimes n} \otimes \cdots\otimes h_{t_k}^{\otimes n}\right) dt_1\dots dt_k
= &\\ \notag
\int \prod_{(i,j) \in J} \rho(t_i - t_j)^{n_{ij}}&1_{[0,T]}(t_i)dt_i 1_{[0,T]}(t_j)dt_j \\
\label{decompo}
\times\int \prod_{(i,j) \notin J} \rho(t_i- t_j)^{n_{ij}}
&\prod_{s=1}^{k-2|J|} 
1_{[0,T]}(t_s)dt_s\,.
 \end{align}
 
 The first factor is actually
  \begin{align}
 \label{prem}\left(\int_{[0,T]^2} \rho(t_2 - t_1)^ndt_1dt_2\right)^{|J|}= \left(\int_{[-T,T]} \rho(t)^n (T-|t|) dt\right)^{|J|}\,,\end{align}
which implies
 \begin{align}
 \label{premlim}
 \lim_{T \to \infty} T^{-|J|}\left(\int_{[0,T]^2} \rho(t_2 - t_1)^ndt_1dt_2\right)^{|J|} = \left(\int \rho(t)^n dt\right)^{|J|}\,.
 \end{align}
 In the second factor, if the complement of $J$ is not empty,
 we denote by 
   $r= k -2|J|$ be the number of remaining blocks.
 
 We can associate to this situation a graph $G$ whose vertices are the $r$ variables of integration and whose edges are the pairs $(k, \ell)$ such that $n_{k\ell} \not=0$. Each vertex has  a degree greater or equal to 2, which implies that the number of edges is greater or equal to $r$. If $G$ is connected, we consider a spanning tree. Such a tree has $r-1$ edges, and each edge can be identified to a linear form
 \[(i,j) \sim (t_1, \dots, t_r) \mapsto t_i-t_j \,.\]
 If these forms are denoted by $\ell_1, \dots, \ell_{r-1}$, the kernel of linear map
 \[(t_1, \dots, t_r) \mapsto (\ell_1, \dots, \ell_{r-1})\]
 is the diagonal $t_1 = \dots = t_r$, hence the rank of the linear map is $r-1$. Adding one coordinate, say  $t_r$, we can perform the change of variables
 \[(t_1, \dots, t_r) \mapsto (u_1= \ell_1, \dots, u_{r-1}= \ell_{r-1}, u_r=t_r)\,.\]
Since its Jacobian is $1$ we get
  \begin{align*}
  \notag
 \int \prod_{(i,j) \notin J} \rho(t_i - t_j)^{n_{ij}}
\prod_{k=1}^{r} 
1_{[0,T]}(t_k)dt_s \leq 
 \int \prod_{(i,j) \notin J} \rho(t_i - t_j)
\prod_{s=1}^{r} 
1_{[0,T]}(t_s)dt_s \\
 \label{compl}
\leq \int_{\mathbb R^{r-1}\times [0,T]}
 \left(\prod_{s=1}^{r-1} \rho(u_s) du_s\right)  du_r 
\leq  
T \left(\int \rho(t) dt \right)^{r-1}\,.
\end{align*}
Now, if $G$ has  $c > 1$ connected components, we cinder each one separately  and we conclude
 \begin{align*}
 \int \prod_{(i,j) \notin J} \rho(t_i - t_j)^{n_{ij}}
\prod_{s=1}^{r} 
1_{[0,T]}(t_k)dt_s \leq T^{c} \left( \int \rho(t) dt\right)^{r -c}\,.
 \end{align*}
Actually , since each component has at least  3 vertices, we have  $3c \leq r$ and we conclude 
 \begin{itemize}
 \item If $|J| = k/2$ then
 \begin{align}
 \label{item1}
\lim_T T^{-k/2}\int_{[0,T]^k} C_\pi \left(h_{t_1}^{\otimes n} \otimes \cdots\otimes h_{t_k}^{\otimes n}\right) dt_1\dots dt_k = \left(\int \rho(t)^n dt\right)^{k/2}
 \end{align}
 \item If $|J| < k/2$, then
 \begin{align}
\notag 
 T^{-k/2}\int_{[0,T]^k} C_\pi \left(h_{t_1}^{\otimes n} \otimes \cdots\otimes h_{t_k}^{\otimes n}\right) dt_1\dots dt_k =&\\
 \label{item2} O(T^{c- k/2+ |J|}) = O(T^{c- r/2}) = O(T^{-r/6})\,.
 \end{align}
 \end{itemize}
 Therefore, the only pair partitions giving a non zero  term as $T \rightarrow  +\infty$ are the inhomogeneous  partitions for which $n_{ij} = n$ for all $i \not=j$, thus a partition of $\Pi_{n,k}$.
 \medskip
 
 3) It is known that a random  symmetric isotropic real Gaussian matrix is of the form $aG + 
  +b \eta I_N $ with 
  \begin{itemize}
  \item $G \el \GOE(N^{-1})$,
  \item 
  $\eta\el \mathcal N(0;1)$  independent of $G$,
  \item $a, b$ real,
  \end{itemize}
  see for instance  \cite[Lemma 4]{chevillard2013} or  \cite[Sec. 2.1]{cheng2018}. The proof extends directly to the Hermitian case with $\GUE$ instead of $\GOE$. Moreover, from the expression of moments (\ref{momentlim}), it is clear that the coefficients $a$ and $b$ are proportional to $\sigma_n$, so that $a_{n,N}$ and $b_{n,N}$ do not depend on $\rho$ and are purely combinatorial.
 
 To characterize $a_{n,N}$ and $b_{n,N}$, it is natural to compute the distribution of $\Tr A M_\infty$ for $A$ a $N\times N$ Hermitian test matrix. On the one hand it is Gaussian with variance
 \begin{align}
 \label{sigma2}
 \sigma_{n,N}^2(A) =\sigma_n^2\left( a^2_{n,N}\tr(A^2)+ b^2_{n,N} (\tr A)^2 \right)\,.
 \end{align}
  On the other hand $\Tr A M_\infty$ as the limit of the stationary Gaussian  process $M_T$ has a variance  
\begin{align}
\label{sigmaint}\sigma_{n,N}^2(A) &:= 2 \int_0^\infty \rho_{A,n,N}(t) dt\\
\label{rhoetr}
   \rho_{A,n,N}(t)&:=
      \E\left(\Tr(AH_\alpha(X_t) \Tr(AH_\alpha(X_0)\right)\,.
 \end{align}
 The following lemma allows to end the proof of Th. \ref{maintheoab} 3.
 \begin{lemma}
 \label{lemrho}
 With $a_{n,N}$ and $b_{n,N}$ defined in (\ref{defa}, \ref{defb}),
 \begin{align}
 \label{eqrho}
 \rho_{A,n,N}(t) = \rho(t)^n \left(a_{n,N}^2 \tr (A^2) + b_{n,N}^2 (\tr A)^2\right)\,.
 \end{align}
 \end{lemma}
 \proof 
 From  definition \eqref{defHermite} of Hermite polynomials as linear combination of tracial monomials of Hermitian Brownian motion, together with formula \eqref{mixmom} for the moments of tracial monomials, we can deduce that $\rho_{A,n,N} (t)$, as a function of $A$,  is a linear combination (depending on $n,N,t$) of terms of the form $\Tr_\beta (D_1, \ldots, D_{2n})$ for some  $\beta \in S_0(2n)$, where all matrices $D_i$ are equal to $\I$ except for two indices where $D=A$. It follows that $\rho_{A,n,N} (t)$ is a linear combination of $(\Tr A)^2$ and $\Tr(A^2)$. Moreover, the coefficients of this combination must be proportional to 
  $\rho(t)^n$.  In the expansion of  $\rho_{A,n,N}(t)$ we have to isolate the terms proportional to $\tr (A^2)\rho^n(t)$ and  $(\tr A)^2 \rho^n(t)$, respectively. 
  
    If we plug  the expansions  (\ref{2var}, \ref{defHermite}) of $H_\alpha(X_t)$ and  $H_\alpha(X_0 )$ in (\ref{rhoetr}) we see that 
  $\rho_{A,n,N}(t)$ is  a linear combination of terms like
  \[\E \left[\Tr \left(A\Tr_{\beta_\pi(\alpha)}(W(h), \dots, W(h))\right) \Tr(A\Tr_{\beta_{\pi'}(\alpha)}(W(h_0)\dots W(h_0)))\right]\]
and using the definition  (\ref{parttrace}) this can be written as
  \[\E\left[\Tr(AW(h)^r) \Tr( A W(h_0)^{s}) \prod_i \Tr\left(W(h)^{v_i}\right)\prod_j \Tr\left(W(h_0)^{v'_j}\right)\right]\,,\]
  for convenient $r,s, v_i, v_j$. Now, owing to (\ref{parttrace}) again, it is of the form
  \begin{align}
  \label{XXAYYA}
  \E\Tr_{\tilde\alpha}\left(X, \dots, X, XA, Y, \dots, Y, YA, X , \dots, X, Y,  \dots, Y\right)\end{align}
    where
    \begin{itemize}
    \item  $X= W(h), Y = W(h_0)$,    
    \item there are $r-1$ successive $X$, then $XA$, then $s-1$ successive $Y$ then $YA$ ...
   \item
    \[\tilde\alpha = (0) (1\dots r)(r+1, \dots, r+s) \gamma \gamma'\]
    with $\gamma$ acting on the last $X$'s and $\gamma'$ acting on the last $Y$'s.
    \end{itemize}
    If we now apply  formula (\ref{mixmom}) we see that the occurrence of $\rho(t)$ arises  from the contractions $C_\pi$. If order to obtain a factor $\rho(t)^n$, the number of variables must be $2n$. \\
    From the definition of $H_{\alpha_n}$, we have $H_{\alpha_n} = X^n + R_n$ where $R_n$ is a trace polynomial of degree strictly less than $n$. Therefore, the only term of the form (\ref{XXAYYA}) containing $2n$ variables (that is $n$ variables $X$, $n$ variables $Y$) is the term 
    $$\E(\Tr(AX^n) \Tr (AY^n)) = \E(\Tr_{\tilde\alpha}\left(X, \dots, X, XA, Y, \dots, Y, YA\right)) $$
     where
    \[\tilde\alpha = (0) (1\dots n)(n+1, \dots, 2n)\,. \]
    Moreover, the term proportional to  $\rho(t)^n$ in  formula (\ref{mixmom}) for $\E(\Tr(AX^n) \Tr (AY^n)) $ is obtained by taking, in the sum, the partitions $\pi \in \mathcal P_2(n,n)$.   For the corresponding term $\Tr_{\pi\tilde\alpha}$ in the RHS of (\ref{mixmom}) we shall use the following auxiliary result. Although straightforward, we state it as a lemma for later reference in Sec \ref{subsub}.
   
   \begin{lemma}
   \label{newlemma}
   For $\tilde\alpha = (0) (1\dots n)(n+1, \dots, 2n)$ and $\pi \in \mathcal P_2(n,n)$ set
   \begin{align*}U_{n,A}^\pi := \mathrm Tr_{\pi\tilde\alpha}[(I_N,\dots, I_N,A, I_N, \dots, I_N, A)
   \end{align*}
    where $A$ is a $N\times N$ Hermitian matrix. Then,
     \begin{itemize}
    \item if  $n$ and $2n$ are in the same cycle of $\pi\tilde\alpha$, then $ U_{n,A}^\pi  = N^{\cyc(\pi\tilde\alpha)} \tr (A^2)\,,$
     \item otherwise 
     \begin{align*}
    U_{n,A}^\pi  = N^{\cyc(\pi\tilde\alpha)} (\tr A)^2\,.    
     \end{align*}
     \end{itemize}
    \end{lemma} 

\subsection{ Asymptotics  $N \rightarrow \infty $}
\label{sec:asymp} 
The limiting distribution $\mu_\infty$ given in Theorem \ref{maintheoab} depends on the index $n$ of the Hermite trace polynomial $H_{\alpha_n}$ and on the dimension $N$ of the matrix. It is natural to keep $n$ fixed in order to recover the $k$-th moment of the $\SC$ distribution. We may then consider either its moments $K_{k,N,n}\sigma_n^k , k \geq 1$ or the coefficients $a_{n,N}, b_{n,N}$ appearing in the decomposition (\ref{second}).

 \subsubsection{The coefficients  $K_{k,n,N}$}
 \label{subsub}
Let us recall that 
$$K_{k, n,N} = \sum_{\pi \in \Pi_{n,k}} C_\pi(\alpha_{nk}) \mbox{ where } C_\pi(\alpha_{nk})  =  \left(\frac{1}{N}\right)^{ nk/2 - \textrm{cyc}_0(\pi \alpha_{nk})}$$.
\begin{proposition}
\label{masterprop}
\begin{enumerate}
\item We have
\begin{align}
\cyc_0(\pi\alpha_{nk}) \leq \frac{nk}{2} 
\end{align}
with equality if and only if the partition $\pi$ in $\Pi_{n,k}$ is non crossing, i.e. the partition on the $n$-blocks of $\pi$ is a non-crossing partition of $[k]$ and satisfies the condition \eqref{bloc-rigide} below. 
\item 
\begin{equation}  \label{Catalan}
\lim_N K_{k, n,N} = C_{k/2}.
\end{equation}
where $C_n$ denotes the $n$th Catalan number.
\end{enumerate}
\end{proposition}
\proof Let $k$ even and $ \pi \in \Pi_{n,k}$. \\
We denote by $\bar 1, \dots, \bar k$ the successive blocks of $[nk]$ ($\bar j = \{(j-1)n +1, \ldots, j n\}$). A partition  $\pi \in \Pi_{n,k}$ induces a pair partition $\bar\pi$ on  the $k$ blocks $\bar j$ in $k/2$ pairs.  Let 
\[\alpha= \alpha_{nk} = (0 1 \dots kn)\]
and denote by  $C_0$ the cycle of  $\pi\alpha_{nk}$ containing $0$. Then, $nk$ belongs to $C_0$. Therefore, the number of elements of $[nk]$, not belonging in $C_0$ is less than $nk-1$, and if 
 $\pi\alpha$ has no singletons,we have :
\[\cyc_0(\pi\alpha_{nk}) < nk/2\,.\]

Assume now that $\pi\alpha_{nk}$  has a singleton. It implies that the partition $\bar\pi$ has an interval $\{\bar j, \overline{j+1} \}$ for some $1 \leq j  \leq k-1$ and $(jn, jn+1) \in \pi$. The singleton of $\pi\alpha_{nk}$ is  $\{jn\}$ and there is no other singletons among the interval $[(j-1)n+1, (j+1)n-1]$. The restriction of $\pi\alpha_{nk}$ on this interval has thus at most $n$ cycles and the maximum of $n$ cycles occurs with a singleton and $n-1$ pairs. \\
If the restriction of $\pi\alpha_{nk}$ on the above interval has a pair $(x,y)$ with $x<y$, it means that $x \leq jn-1$, $y \geq  jn+2$ and $\pi(x+1) = y, \ \pi(y+1) = x$. Thus, the only way for the restriction of $\pi$ to have $n$ cycles (a singleton and $(n-1)$ pairs) is that $\pi$ is non-crossing on this interval, i.e. satisifies :
$$   \pi((j-1)n + r) = (j+1)n +1 -r, \, 1 \leq r \leq n; $$
We now prove that $(j-1)n$ and $(j+1)n$ are in the same cycle $C_j$ of $\pi \alpha_{nk}$. Note that $ a:= (\pi\alpha)^{(-1)}((j-1)n)$  does not belong to $(\bar{j} \cup \overline{j+1})$. The image $\pi\alpha((j-1)n) = \pi((j-1)n+1) \in (\bar{j} \cup \overline{j+1})$. The only way to pass from $\pi\alpha((j-1)n)$ to $a$ in the cycle $C_j$ (i.e. to exit $\bar{j} \cup \overline{j+1}$) is to pass by the point $(j+1)n$.
We can thus "identify" $(j-1)n$ and $(j+1)n$ to count the number of cycles. Thus, we are led to the study of the restriction of $\pi$ on $[ kn] \backslash (\bar{j} \cup \overline{j+1})$ with $\tilde{\alpha_{nk}} = (012\ldots (j-1)n \ (j+1)n+1 \ldots kn)$. A recursive argument on $k$  gives $\cyc_0(\pi\alpha_{nk}) \leq nk/2\,$.
Therefore, the number of partitions $\pi$ giving the maximum number of cycles $nk/2$ is exactly the number of non-crossing pair partitions $\bar\pi$  (recall than a non-crossing partition contains at least one interval $(\bar{j}, \overline{j+1})$) which satisfy the following condition :
if $(\bar j, \bar l)$, $j <l$ is a pair of $\bar \pi$, then,
\begin{equation} \label{bloc-rigide}
 \pi((j-1)n + r) = ln +1 -r, \, 1 \leq r \leq n.
 \end{equation}
For such a partition, $\lim_N C_\pi(\alpha_{nk}) = 1$ and $\lim_N C_\pi(\alpha_{nk}) = 0$ in the other cases.
The number of partitions $\pi$ in $\Pi_{n,k}$ giving the maximal number of cycles is thus equal to the number of non-crossing pair partitions $\bar \pi$, i.e. the Catalan number $C_{k/2}$, proving \eqref{Catalan}. 
$\Box$

\subsubsection{The coefficents $a_{n,N}$ and $b_{n,N}$}

\begin{proposition}
\label{prop:abN}
\begin{enumerate}
\item
The coefficients $a_{n,N}$ and $b_{n,N}$ satisfy :
\begin{align}a^2_{n,N} = 1 + O(N^{-1})\ , \quad b^2_{n,N} =n-1 + O(N^{-1})\,.\end{align}
\item As $N\to \infty$, the empirical spectral distribution of $M_\infty$ converges in distribution to $\SC(0 ; \sigma^2_n)$.
\end{enumerate}
\end{proposition}
The first statement will be a consequence of Lemma \ref{termedominant} below. For the second one, it is enough to observe that for fixed $n$,
\begin{itemize}
\item
the empirical spectral distribution of $\sigma_n b_{n,N}\xi I_N$ converges to $\delta_0$  as $N\to \infty$, 
\item
the empirical spectral distribution of  $\sigma_n a_{n,N}\GUE(N^{-1})$ converges to $\SC(0; \sigma^2_n)$  as $N\to \infty$,
\item $\sigma_n b_{n,N}\xi I_N$ and  $\sigma_n a_{n,N}\GUE(N^{-1})$ are independent.
\end{itemize}

\noindent {\bf Proof of Proposition \ref{prop:abN}} 1
\smallskip

With the notations of Lemma \ref{newlemma}, we see that if $\pi \in \mathcal P_2(n,n)$, then $\pi\tilde\alpha$ has no singleton. Therefore 
$U_{n,A}^\pi = O(N^n)$
and a contribution of order $N^n$ arises precisely  when $\pi\tilde{\alpha}$ is a product of $n$ pairs ($\times (0)$). \\
We  study those pairings $\pi$ such that $\pi\tilde{\alpha} = (0) \prod_{i=1}^n  \tau_i$ where every $\tau_i$ is a pair (or transposition). Two cases appear : \\
1) one  $\tau_i$ is $(n,2n)$. In this case, 
$$ U_{n,A}^\pi
 = N^{n-1} \Tr(A^2) = N^n \tr(A^2). $$
2) $n$ and $2n$ are not in the same pair of $\pi$. Then,
$$ U_{n,A}^\pi
 = N^{n-2} (\Tr A)^2 = N^n (\tr A)^2. $$
It remains to compute the number of  such pairings $\pi$.

\begin{lemma} \label{termedominant}
1) There exists a unique partition $\pi$, namely the non-crossing partition $\pi = (1, 2n) (2, 2n-1) \ldots (n, n+1)$ of $ \mathcal P_2(n,n)$, such that
$$ 
U_{n,A}^\pi  = N^n \tr(A^2).$$
2) There are $n-1$ pair partitions $\pi \in \mathcal P_2(n,n)$ such that :
$$
U_{n,A}^\pi  = N^n (\tr A)^2. $$
\end{lemma}

\noindent {\bf Proof of Lemma \ref{termedominant}} 
\begin{enumerate}
\item We assume that $\pi\tilde{\alpha} = (1,j_1) (2,j_2) \ldots (n-1, j_{n-1})(n,2n)$. Note that 
$$\tilde{\alpha}^{-1} = (0) (n \ (n-1) \ldots 1)( 2n \ (2n-1) \ldots n+1).$$
We have successively
\begin{align*}
\notag
\pi(1) =\pi\tilde{\alpha}\tilde{\alpha}^{-1} (1) = \pi\tilde{\alpha}(n) = 2n\\
\pi(2n)=\pi\tilde{\alpha}\tilde{\alpha}^{-1} (2n) = \pi\tilde{\alpha}(2n-1)\,,
\end{align*}
but since $\pi$ is a pairing $\pi(2n) =1$ so that
\begin{align*}
\pi\tilde{\alpha}(2n-1) = 1\,,
\end{align*}
which implies $j_1 = 2n-1$.
Then, similarly 
\begin{align*}
\notag
\pi(2) &= \pi\tilde{\alpha}\tilde{\alpha}^{-1}(2)= \pi\tilde\alpha(1) = j_1 = 2n-1\\
\notag
2= \pi(2n-1) &=\pi\tilde{\alpha}\tilde{\alpha}^{-1}(2n-1) = \pi\tilde{\alpha}(2n-2)\\
j_2&= 2n-2
\end{align*}
and so on. Eventually
\begin{align*}
\pi= (1,2n)(2, 2n-1)\dots(n, n+1)\,
\end{align*}
so that there is one and only one $\pi$ giving the dominant term of $\tr(A^2)$.

\item We now assume that 
 \[\pi\tilde{\alpha}= (i_1, j_1) \dots (i_{n-2}, j_{n-2})(n, j_{n-1}) (i_n, 2n)\]
We set 
\begin{align}
\notag
\pi(n) &= \pi\tilde{\alpha}(n-1)=: k \in [n+1, 2n]
\end{align}
If $k= n+1$ we have 
\begin{align*}
\pi(n) &= n+1 = \pi\tilde{\alpha}\tilde{\alpha}^{-1}(n) = \pi\tilde{\alpha}(n-1)\\
\pi(n+1) &= n = \pi\tilde{\alpha}\tilde{\alpha}^{-1}(n+1) = \pi\tilde{\alpha}(2n)
\end{align*}
This is a contradiction since, by assumption, $(n, 2n)$ is not a pair of $ \pi\tilde{\alpha}$. \\

Therefore,   $k= n+1+j$ for $j\in [1, n-1]$ and we have successively:
\begin{align*}
\pi(n) &= n+1+j = \pi\tilde{\alpha}\tilde{\alpha}^{-1}(n) = \pi\tilde{\alpha}(n-1)\\
\pi(n+j+1) &= n = \pi\tilde{\alpha}\tilde{\alpha}^{-1}(n+j+1) = \pi\tilde{\alpha}(n+j)
\end{align*}
Thus $(n, n+j) \in \pi\tilde{\alpha}$ and 
$$ \pi(1)  = \pi\tilde{\alpha}\tilde{\alpha}^{-1}(1) = \pi\tilde{\alpha}(n)= n+j.$$
implying that $1 = \pi(n+j) = \pi \tilde{\alpha} (n+j-1)$. \\
Thus $\pi(2) = \pi \tilde{\alpha}(1) = n+j-1$. 
We obtain successively:
\begin{align*}
\pi(l) &= \pi\tilde{\alpha}\tilde{\alpha}^{-1}(l-1) = n+j -l +1, \; l\leq j \\
&\cdots\\
\pi(j) &= \pi\tilde{\alpha}\tilde{\alpha}^{-1}(j-1) = n+1\\
\pi(n+1) &= j = \pi\tilde{\alpha}\tilde{\alpha}^{-1}(n+1) = \pi\tilde{\alpha}(2n)\\
\pi(j+1)  &= \pi\tilde{\alpha}\tilde{\alpha}^{-1}(j+1) = \pi\tilde{\alpha}(j)=2n\\
\pi(2n) &= j+1 = \pi\tilde{\alpha}\tilde{\alpha}^{-1}(2n) = \pi\tilde{\alpha}(2n-1)\\
&\cdots\\
\pi(n-1) &= \pi\tilde{\alpha}\tilde{\alpha}^{-1}(n-2) = n+j+2\,.
\end{align*}
and eventually
\[\pi = (1, n+j) (2, n+j-1) \dots (j, n+1) (j+1, 2n)(j+2, 2n-1)\dots(n, n+j+1)\,.\]
In conclusion, there are exactly $n-1$ pairings giving the dominant term in $(\tr A)^2$.
\end{enumerate}
\begin{remark}
One  may also ask what happens when $N\to \infty$ for the expression of Hermite trace polynomials studied above. For the four examples given in Proposition \ref{4H}, if we admit that the terms involving $\tr$ converge to the moments of the $\SC$ distribution, we obtain -up to a sign change- the first four Chebyshev polynomials. 
This also holds for the fifth one, although we omit its explicit expression.
\end{remark}

\section{Matrix-variate Hermite polynomials}
\label{mvariate}
\subsection{Introduction}
In this section, we  consider two models of matrix-valued processes $X(t)$ obtained by smoothing a matrix Brownian motion. As in the previous sections, we 
study the asymptotic behavior of  $\ve^{1/2}\int_0^{1/\ve} F(X(t)) dt$, but now $F$ will be  a suitable real function such that $F(X(t))$ is centered. The classical key tool is the decomposition according to some Hermite-like basis, called the matrix-variate Hermite polynomials.

In the first model, we consider the bilateral Brownian motion $\mathcal W\sn(t)$ on the space $S_N$ of $N\times N$ symmetric real matrices \footnote{These Brownian motions are not normalized.}\edef\fn{\the\value{footnote}} and a smoothing real function $\varphi$ satisfying $\Vert \varphi\Vert_2 = 1$.

 In the second model, we consider the bilateral Brownian motion $\mathcal W^{\ell, N}(t)$ on the space $\R^{\ell\times N}$ of $\ell \times N$ rectangular real matrices\footnotemark[\fn]and a smoothing function $\Phi$ living in   the space of $N\times N$ matrices and  satisfying
 \begin{align}
 \label{crucialrectangular}
 \int \Phi(s)^T\Phi(s) ds = \hbox{Id}_N\,.
 \end{align}

 The smoothed processes are defined respectively by
\begin{align}
\label{basegen}
X(t) &= \int_{-\infty}^\infty \varphi(t-s) d\mathcal W\sn (s) \in S_N\,,\\
\label{rectgen}
X(t) &= \int_{-\infty}^\infty  d\W^{(\ell,N)} (s)\Phi(t-s) \in \R^{\ell\times N}\,.
\end{align}
We use the same notation $X(t)$ for both models, the context making clear which one is meant. 
These processes are stationary Gaussian 
 and  their covariances are given respectively by
 \begin{align}
 \label{covs}
 \E \left( X(t_2)_{ij}X(t_1)_{rs}\right) = \rho(t_2 -t_1) \left(\delta_{ir}\delta_{js}+ \delta_{is}\delta_{jr}\right)\,,
 \end{align}
where $\rho$ was defined in (\ref{defcov}) and
\begin{align}
 \label{covr}
 \E\left (X(t_2)_{ij}X(t_1)_{rs}\right) = \delta_{ir} R(t_2 - t_1)_{js}
 \end{align}
 where
 \begin{align}
 \label{defR}
 R(t) =   \int\Phi(s)^T \Phi(t+s) ds \in \R^{N\times N}\,.
 \end{align}
 The definition of  matrix-variate Hermite polynomials in both the symmetric and rectangular cases, relies on zonal polynomials, which we briefly recall  below (see also \cite{chikuse92}). 
  
  Let $Z$ be a symmetric $m\times m$ matrix, let 
   $n$ be an integer and $\kappa$ be a partition of $n$  in no more than $m$ parts denoted
 $\kappa \vdash n$, i.e. $\kappa = (k_1, \dots, k_l),\  l \leq m,\  k_1 \geq k_2\geq k_l\geq 1$, and $n= k_1 + \dots+ k_l$.   
 The zonal polynomials $C_\kappa(Z)$ are defined as a basis of the space of all homogeneous symmetric polynomials in the eigenvalues of $Z$.  
 
 Actually, to simplify the exposition, the main results, in Propositions \ref{prop:12} and \ref{proprect} are stated only for polynomials. 

\subsection{Symmetric real matrices}
The Hermite polynomials $\He_\kappa\sn (X), \kappa\vdash n =0, 1, \cdots$, where $\kappa$ is a partition with no more than $N$ parts, constitute a complete system of orthogonal polynomials associated with the distribution of density
\begin{align}
\label{GOE}
\phi_N(X) = 2^{-N/2}\pi^{-N(N+1)/4}\exp \left(- \frac{1}{2}\Tr(X^2) \right)
\end{align}
with respect to the Lebesgue measure on $\R^{N(N+1)/2}$.
Actually
\begin{align*}
\int_{S_N} \He_\kappa\sn (X) \He_\sigma\sn (X) \phi_N (X) \dd X = \delta_{\kappa, \sigma} n ! C_\kappa(\I)
\end{align*}
(\cite[formula (2.13)]{lawi2008}, \cite[formula 5.3]{chikuse2020}).

We assume that  $\mathcal W\sn(t)$ is a real symmetric matrix-valued Brownian  process, so that $X(t)$ defined in (\ref{basegen}) is a stationary process with 1d-marginal density $\phi_N$.  We set
\begin{align}\zeta_\ve^{(\kappa)} &= \ve^{1/2}\int_0^{\ve^{-1}} \He\sn_\kappa(X(t)) dt\,.\end{align}
The main result of this section is the following proposition.
\begin{proposition} 
\label{prop:12}
As $\ve \to 0$,
\begin{align}
\zeta_\ve^{(\kappa)} \vers_{\ve \to 0}^{(d)} \mathcal N(0; \sigma_{N, \kappa}^2)
\end{align}
where
\[\sigma_{N, \kappa}^2 = 2\ n ! C_\kappa(\I) \int_0^\infty \rho(t)^n dt\,,\]
and $\rho$ was defined in a previous section.
\end{proposition}
 \proof  
We skip the upper index $\sn$. The centered process $\left(\He_\kappa(X(t)), t\geq 0\right)$ is a non-linear function  of a Gaussian multi-dimensional process and we apply the  Breuer-Major theorem.
We will prove that the covariance of the process $\He_\kappa(t)$ defined by
\[\Gamma(t) = \mathbb E\left(\He_\kappa(X(t))\He_\kappa(X(0)\right)\]
satisfies
\begin{align}
\label{611}
\Gamma(t) = \rho(t)^n n ! C_\kappa(\I)\,.
\end{align}
If $t$ is fixed, 
we observe the identity in distribution
\begin{align}
\label{idinlaw}
\left(X(t), X(0)\right) \el \left(\rho(t) X(0) + \sqrt{1 - \rho(t)^2}X', X(0)\right)
\end{align}
where $X(0)$ and $X'$ are independent and have density $\phi_N$.

Summarizing, we have to compute
 $$\mathbb E\left( \He_\kappa(X)\He_\kappa(Y)\right)$$ with
\begin{align}
\label{Y=}
Y = \rho X + \sqrt{1-\rho^2} X'\ , \ X\perp X', X\el X'\,.\end{align}
We will condition upon $X$ and use a reduction to a computation coordinatewise.
We will use the following lemma, which is probably well known.
\begin{lemma}There exists a system of numerical constants $u_\alpha^\kappa$ indexed by  $\kappa$, a partition of $n$ and   $\alpha= (\alpha_{ij}\in \mathbb N^{\frac{N(N+1]}{2}})$,  
a multi-index with  $|\alpha| =\sum_{i\leq j} \alpha_{ij} = n$ such that
\begin{align}
\label{prod1}
\He_\kappa(X) = \sum_{|\alpha|= n}u_\alpha^\kappa 
\left(\prod_i  2^{-\alpha_{ii}/2} H_{\alpha_{ii}}(X_{ii}\sqrt 2)\right) \prod_{i< j} H_{\alpha_{ij}}(X_{ij})\,,
\end{align}
where the polynomials $H_p$, $p\geq 1$ are the classical (scalar) Hermite polynomials.
\end{lemma}
\proof 
 
The proof is based on the Rodrigues formula  in the symmetric case is \cite[Sec. 3.2]{chikuse92}. It says that if  
\begin{align*}
\partial X =: \left(\frac{1}{2}(1 + \delta_{ij}) \frac{\partial}{\partial x_{ij}}\right)_{i,j}\,.
\end{align*}
then 
\begin{align}
\He_\kappa (X) = \phi_N(X)^{-1} C_\kappa(-\partial X) \phi_N (X)\,.
\end{align}

Besides we know that,  if $Z \in S_N$ then  by definition, $C_\kappa (Z)$ is a symmetric  homogeneous polynomial in the eigenvalues. It has then a decomposition
\[C_\kappa (Z) = \sum_\nu
 c_{\kappa \nu}\ s_1(Z) ^{\nu_1}\dots s_k(Z)^{\nu_k}\,,\]
where $s_ j(Z) = \tr Z^j$ and $\nu = (1^{\nu_1} 2^{\nu_2}\dots k^{\nu_k})$
with $\sum_j  j \nu_j =  n$. Since $t_j(Z)^{\nu_j}$ is a homogeneous polynomial of degree $j\nu_j$ in the entries of $Z$, we conclude that $C_\kappa (Z)$ has  
a decomposition
\begin{align} 
\label{decomp}
C_\kappa(Z) = \sum_{|\alpha|= n}u_\alpha^\kappa \prod_{i\leq j} Z_{ij}^{\alpha_{ij}}\,,
\end{align}
(see also \cite[formula 2.6]{notar2023}).

This leads to
\begin{align*}
\He_\kappa (X) = \phi_N (X)^{-1}
\sum_{|\alpha|= n}(-1)^n
u_\alpha^\kappa  &\prod_{i< j}2^{-\alpha_{ij}}
\frac{\partial^{\alpha_{ij}}}{\partial X_{ij}^{\alpha_{ij}}}
&\prod_i \frac{\partial^{\alpha_{ii}}}{\partial X_{ii}^{\alpha_{ii}}}\phi_N(X)\,.
\end{align*}
Now, the obvious decomposition
\[\phi_N(X) = \prod_i \phi_1(X_{ii})\prod_{i<j} \phi_1(X_{ij}\sqrt 2)\,,\]
gives
\begin{align}
\notag
\He_\kappa (X) =  \sum_{|\alpha|= n}(-1)^ n
u_\alpha^\kappa
 &\prod_{i<j}2^{-\alpha_{ij}}\phi_1(X_{ij}\sqrt 2)^{-1}\frac{\partial^{\alpha_{ij}}}{\partial X_{ij}^{\alpha_{ij}}}\phi_1(X_{ij} \sqrt 2)\\
&\times \prod_i \phi_1(X_{ii})^{-1}\frac{\partial^{\alpha_{ii}}}{\partial X_{ii}^{\alpha_{ii}}}\phi_1(X_{ii})
\end{align}
Then, using the classical definition of 
classical Hermite polynomials :
\[\phi_1(x)^{-1}\frac{d^p}{dx^p}\phi_1(x) = (-1)^p H_p(x)\,,\]
we get
 (\ref{prod1}).
\bigskip

\noindent\underline{End of the proof of Prop. \ref{prop:12}}

In the setting of (\ref{Y=})
 and using the property of the Ornstein-Uhlenbeck semigroup, we get
\begin{align}
\mathbb E\left[H_{\alpha_{ij}}(\rho X_{ij} + \sqrt{1-\rho^2}X'_{ij}) | X\right]=\rho^{\alpha_{ij}}H_{\alpha_{ij}}(X_{ij})\,.
\end{align}
A plugging into (\ref{prod1}) leads to 
\begin{align}
\notag
\mathbb E [\He_\kappa (Y) | X] &= \sum_{|\alpha|= n}u_\alpha^\kappa  \left(\prod_i  2^{-\alpha_{ii}/2} \rho^{\alpha_{ii}}H_{\alpha_{ii}}(X_{ii})\right)\prod_{i< j}\rho^{\alpha_{ij}} H_{\alpha_{ij}}(X_{ij})\\
\notag
&= \rho^n \sum_{|\alpha|= n}u_\alpha^\kappa \left(\prod_i  2^{-\alpha_{ii}/2} H_{\alpha_{ii}}(X_{ii})\right) \prod_{i<j} H_{\alpha_{ij}}(X_{ij})
\\
\notag
&= \rho^n \He_\kappa (X)\,,
\end{align}
and this proves that (\ref{611}) holds true.
\begin{remark}
Another proof of the Lemma could be based on  a striking representation \cite[(5.48) p. 100]{chikuse2020} :
\begin{align}
\label{careful}
\He_\kappa(X) =
\int _{S_N}  C_\kappa(X+i M)  \phi_N(M) \dd M\,.
\end{align}
  The second ingredient would be the decomposition (\ref{decomp}) 
and independence of coordinates. The issue is that we don't know any definition of $C_{\kappa}$ applied to a complex (although symmetric) matrix.
\end{remark}

\subsection{Hermitian matrices}
A similar analysis could be made for Hermitian matrices as well. We do not give details, but present some examples.

The following proposition, taken from  \cite[Cor. 6.17]{anshelevich2022}, gives 
 an explicit way to compute the matrix-variate Hermite polynomials in the Hermitian context. Following the notations therein, if  $\kappa$ is a partition of $n$, let us  denote by $\chi^\kappa$ the character of the irreducible representation of $S(n)$ corresponding to $\kappa$ and let us identify $\chi^\kappa$ with the element
\[\sum_{\sigma \in S(n)} \chi^\kappa(\sigma) \sigma  \in \mathbb C[S(n)]\,.\]
\begin{proposition}\cite[Cor. 6.17]{anshelevich2022}
The Hermite polynomial of matrix argument $X$, indexed by a partition $\kappa$, is a multiple of
\begin{align}
\label{additional}
\widetilde{{\mathbf H}}_\kappa\sn (X)  := \sum_{\sigma \in S(n)} \chi^\kappa (\sigma)\tilde{H}_{\hat\sigma} (X, N)
\end{align}
where $\tilde{H}_\alpha (X, u)$ is defined by \eqref{2var}\footnote{In \cite[Cor. 6.17]{anshelevich2022}, the above condition $u= N$ was improperly replaced by $u =1/N$.} and $\hat\sigma \in S_0(n)$ is the permutation $\sigma$ with the additional cycle $(0)$.
\end{proposition}

For the first values $n$ we state the following proposition, whose proof is in Section \ref{sec:app}.
\begin{proposition}
\label{prop:sigmahat}
\begin{align}
\label{621}
\hat\sigma= (0)(12)\ &,  \ \tilde{H}_{\hat\sigma} (X, N) =\Tr (X^2) -N^2\\
\label{622}
\hat\sigma= (0)(1)(2) \ &, \ \tilde{H}_{\hat\sigma} (X, N) =(\Tr X)^2 - N\\ 
\label{123}
\hat\sigma = (0)(123) \ &, \ \tilde{H}_{\hat\sigma} (X, N) = -\Tr (X^3)+ 3N\Tr X\\
\label{12}
\hat\sigma = (0)(12)(3) \ &, \ \tilde{H}_{\hat\sigma} (X, N) = (N^2+ 2)\Tr X - \Tr (X^2) (\Tr X)\\
\label{id}
\hat\sigma= (0)(1)(2)(3) \ &, \ \tilde{H}_{\hat\sigma} (X, N) = 3 N\Tr X - (\Tr X)^3\,.
\end{align}
\end{proposition}
The coefficients $\chi^\kappa (\sigma)$ associated to the partitions  (or Young tableaux)
  are given in  particular in  
 \cite[Ex. 6 p. 14 and Ex. 4.5 p. 47]{fulton2013}.
 
1)  If $\kappa=(1,1)$ then $\chi^\kappa = \mathrm{id} - (12)$ so that  from (\ref{additional}) 
 and (\ref{621} \ref{622})
 
 \begin{align}
 \label{H11}\widetilde{\bf H}_\kappa\sn (X) = \left((\Tr X)^2 - N\right) -\left(\Tr (X^2) - N^2\right)  \,.\end{align}
 
 2) If $\kappa = (2)$ then $\chi^\kappa = \mathrm{id} + (12)$ so that
 \begin{align}
 \label{H2}
 \widetilde{\bf H}_\kappa\sn (X) = \left((\Tr X)^2 - N\right) +\left(\Tr (X^2) - N^2\right) \,.\end{align}
 We can check easily  that these polynomials are orthogonal under $\phi_N$.

3)  If
$\kappa =(1,1,1)$ then
$\chi^\kappa = \hbox{id} - (12) - (13)-(23) +(123) + (132)$ so that  from (\ref{additional}) and  (\ref{123}-\ref{12}-\ref{id})
\begin{align*}
\widetilde{{\mathbf H}}_\kappa\sn (X)
= -(\Tr X)^3  + 3 (\Tr X) \Tr (X^2) - 2 \Tr (X^3) 
 -3(N-1)(N-2)\Tr X\,.
\label{H111}
\end{align*}

4) If 
$\kappa= (2,1)]$, then $\chi^\kappa =2 \hbox{id} -(123) - (132)$ so that 
\begin{align*}\label{H21} 
\widetilde{{\mathbf H}}_\kappa\sn (X) &= 2\Tr (X^3) -  2 (\Tr X)^3\,.
\end{align*}

5) If
$\kappa = (3)$  then $\chi^\kappa= \hbox{id} + (12) + (13)+(23) +(123) + (132)$ so that 
\begin{align*}
\label{H3}
\widetilde{{\mathbf H}}_\kappa\sn (X) & =-(\Tr X)^3 -3( \Tr X)\Tr( X^2) - 2 \Tr (X^3) + 3(N+1)(N+2) \Tr X\,..
\end{align*}

These results are consistent with the expressions of Hermite polynomials given in \cite{dumitriu2007}. Therein,  the complex case corresponds to $\alpha = 2/\beta =1$, $C^1_\kappa$ are the complex zonal polynomials and $s_i := \Tr(X^i)$.
\begin{align*}
\notag
H^{1}_{(1,1)}&=  C^1_{(1,1)} + \frac{N(N-1)}{2} \\
&= \frac{1}{2}\left(s_1^2 - s_2+ N(N-1)\right)\,.\\
\notag
H^{1}_{(2)}&= C^1_{(2)} - N(N+1)\\
&= s_1^2 + s_2 -N(N+1)\,.\\
\notag
H^{1}_{(1,1,1)} &=  C^1_{(1,1,1)} + \frac{(N-1)(N-2)}{2} C_{(1)}^1\\
 &= \frac{1}{6}s_1^3 - \frac{1}{2}s_1s_2 + \frac{1}{3} s_3+ \frac{(N-1)(N-2)}{2} s_1\,.\\
\notag
H^{1}_{(2,1)}&= C^1_{(2,1)}  + 0 \  C_{(1)}^1\\
&=  \frac{1}{3}\left(s_1^3 - s_3\right)\,.\\
\notag
H^{1}_{(3)}&=  C_{(3)}^1 + \frac{(N+1)(N+2)}{2} C_{(1)}^1\\
& = \frac{1}{6} s_1^3 + \frac{1}{2}s_1s_2 + \frac{1}{3}s_3 + \frac{(N+1)(N+2)}{2}s_1\,.
\end{align*}

\subsection{Rectangular matrices}

We assume in this part that $\ell\leq N$. The Hermite polynomials  $\He_\kappa^{(\ell,N)}, \kappa \vdash n =0,1, \dots$, where $\kappa$ is a partition in no more than $\ell$ parts, and whose variates are rectangular $\ell\times N$ real matrices constitute a complete system of orthogonal polynomials associated with the distribution $\mathcal N_{\ell\times N}(0, \mathrm{Id}_\ell \otimes \I)$ of density 
\begin{align}
\label{rect}
\phi_{\ell, N}(X) = (2\pi)^{-\ell N/2} \exp\left(-\frac{1}{2}\Tr (XX^T)\right)
\end{align}
with respect to the Lebesgue measure on $\R^{\ell N}$. Actually
\begin{align}
\int_{\R^{\ell N}}\He_\kappa^{(\ell,N)}(X) \He_\sigma^{(\ell,N)}(X)\phi_{\ell, N}(X) = \delta_{\kappa, \sigma } 4^n \left(\frac{N}{2}\right)_\kappa n ! C_\kappa(\mathrm{Id}_\ell)\,. 
\end{align}
We refer to \cite[Sec. 4]{chikuse92} and \cite{notar2023} for the notations and for the interest of this model.  

We assume that  $\mathcal W^{(\ell,N)}(t)$ is a  $\ell\times N$ matrix-valued  process whose coordinates  are independent bilateral standard Brownian motions.  
Notice that owing to the assumptions (\ref{crucialrectangular}) and (\ref{covr}),  the  process $X(t)$ defined in (\ref{rectgen}) has marginals distributed as  $\N_{\ell\times N}(0, \mathrm{Id}_\ell\otimes \I)$.

Recall the definition (\ref{crucialrectangular})
\[R(t) = \int \Phi(s)^T\Phi(t+s) ds\] 
which satisfies by assumption (\ref{crucialrectangular}) $R(0)  = {\mathrm Id}_N$.
In this  context  the correlation function $\rho$ will be replaced by the symmetric, semidefinite positive matrix
\begin{align}
{\mathbf R}(t)  = (R(t)R(t)^T)^{1/2}\,.
\end{align}
 Actually, we have like {(\ref{CSrho})
 \begin{align}
 \label{unitball}
 {\mathbf R}(t) \leq {\mathrm Id}_N\,.
 \end{align}
It is a consequence of the matrix Cauchy-Schwarz inequality - see for instance \cite[formula (2.9)]{albadawi}- which  says that if $|A| := (A^T A)^{1/2}$ then 
\[\left|\sum_1^m B_i A_i \right|\leq\left(\sum_1^m |A_i|^2\right)\]
whenever $\sum_1^m |B_i^T|^2 \leq {\mathbf 1}$. 
This inequality can be easily extended to integrals instead of sums, and taking  $A_i \equiv \Phi(s)^T$, $B_i \equiv \Phi(t+s)$
we get exactly $|R(t)| \leq \mathrm{Id}_N$ hence  (\ref{unitball}) since $R(t) R(t)^T$ and $R(t)^TR(t)$ have the same nonzero eigenvalues.

We set
\begin{align}
 \eta_\ve^{(\kappa)} =\ve^{1/2}\int_0^{1/\ve} \He_\kappa^{(\ell, N)}(X(t))dt\,,
\end{align}
The main result of this section is the following proposition.

\begin{proposition}
\label{proprect}
As $\ve\to 0$
\begin{align}
\eta_\ve^{\kappa} \vers_{\ve \to 0}^{(d)} \mathcal N(0; \sigma_{\ell, N, \kappa}^2)
\end{align}
where
\[\sigma_{\ell, N, \kappa}^2 = 2 \times  4^{- n}  n ! \left(\frac{N}{2}\right)_\kappa \frac{C_\kappa(\mathrm{Id}_\ell)}{C_\kappa (\I)} \int_0^\infty C_\kappa({\mathbf R}^2(t))  \dd t\,.\]
\end{proposition}
\proof
In the sequel, we skip the index $^{(\ell, N)}$ for simplicity. 

Since  $X(t)$ is a standard Gaussian matrix and since $\He_0 \equiv 1$, the process  
$\left(\He_\kappa(X(t)), t\geq 0\right)$ is centered.  Actually it is is a non-linear function of a Gaussian multi-dimensional real process. We apply a continuous version of the 
 Breuer-Major theorem. 
We have to prove that the  covariance of  $\He_\kappa(t)$ 
satisfies
\begin{align}
\label{112}
\mathbb E\left(\He_\kappa(X(t))\He_\kappa(X(0)\right) =  4^{- n} n ! \left(\frac{N}{2}\right)_\kappa \frac{C_\kappa(\hbox{Id}_\ell)}{C_\kappa (\I)}  C_\kappa({\mathbf R}^2(t))\,.
\end{align}
We continue as in the end of proof of Prop. \ref{prop:12}, with the help of \cite[Th. 3.12]{notar2023}, but with a careful treatment of $R(t)$ which is not necessarily symmetric. 
\begin{lemma}
Set ${\mathbf S}(t)= \left(\I- {\mathbf R}(t)^2\right)^{1/2}$. Then for fixed $t$, 
\begin{align}
\left(X(t), X(0)\right) \el \left(Y, X(0)\right)\,,
\end{align}
where $Y:=X(0)R(t)^T+ X'{\mathbf S}(t)$  with $X'\el X(0)$
 and independent of $X(0)$.
\end{lemma}
\proof

We have, by independence 
\begin{align}
\notag
\E(Y_{ij}Y_{rs}) &= \E\left((X(0)R(t)^T)_{ij}(X(0)R(t)^T)_{rs}\right) + \E\left((X'S(t))_{ij}(X'S(t))_{rs}\right)\\
\notag
&= \sum_{k,p} \left(R(t)_{jk}R(t)_{sp} + {\mathbf S}(t)_{jk}{\mathbf S}(t)_{sp}\right)\delta_{ir}\delta_{kp} 
\\
\notag
&
= \left(R(t)R(t)^T + {\mathbf S}(t)^2 \right)_{js}\delta_{ir} 
= \left({\mathbf R}(t)^2 + {\mathbf S}(t)^2 \right)_{ir}\delta_{js} = \delta_{ir}\delta_{js}\,,
\end{align}
which are the covariances of $X(t)$. 

Moreover, again by independence 
\begin{align*}\notag
\E\left(Y_{ij}X(0)_{rs}\right)&= \E\left((X(0)R(t)^T)_{ij}X(0)_{rs}\right)\\
&= \sum_k R(t)_{jk}\delta_{ir}\delta_{ks}= R(t)_{js}\delta_{ir}\,,
\end{align*}
which are the cross covariances of $X(t)$ and $X(0)$, as given in (\ref{covr}). $\Box$
\bigskip

\noindent\underline{End of the proof of Prop. \ref{proprect}}

From the above lemma, we have
\begin{align*}
\mathbb E\left(\He_\kappa(X(t))\He_\kappa(X(0)\right) = \mathbb E\left(\He_\kappa(Y)\He_\kappa(X(0)\right)
\end{align*}
but, since for any $N\times N$ deterministic orthogonal matrix $X(0) \el X(0)U$, we get
\begin{align}
\mathbb E\left(\He_\kappa(X(t))\He_\kappa(X(0)\right) = \mathbb E\left(\He_\kappa(X(0)UR(t)^T + X' {\mathbf S}(t))\He_\kappa(X(0)U\right)\,. 
\end{align}
Now, it is known that $\He_\kappa(MU) = \He_\kappa(U)$  for any $M \in \R^{\ell\times N}$, so that choosing $U= \left( {\mathbf R}(t)\right)^{-1} R(t)$ (polar decomposition of $R(t)$), 
we conclude that
\begin{align*}
\mathbb E\left(\He_\kappa(X(t))\He_\kappa(X(0)\right) = \E\left(\He_\kappa(X(0){\mathbf R}(t)^T + X' {\mathbf S}(t))\He_\kappa(X(0)\right)\,. 
\end{align*}
It is then enough to apply 
  \cite[Th. 3.12]{notar2023}. $\Box$
 
\section{Appendix : Explicit computations of Sections \ref{sec:matrix2} and \ref{mvariate}}
\label{sec:app}
In this section we
 refer to definition \ref{defiAns} and   set $q= 1/N$, $c= a(\alpha, \pi)$ and $\beta = P^{[0,n] \backslash \textrm{supp}(\pi)} _{[0,n-2l]} (\pi \alpha) \big|_{\textrm{supp}(\pi)^c}$.
\medskip

\subsection{Proof of Proposition \ref{4H}}
\label{appendix}

We apply formulas (\ref{2var}) and (\eqref{defHermite}).
 
\noindent For $n =1$, $\mathcal P_{1,2}(n) = (1),l = 0,  n-|\pi|= 1, \pi\alpha = (01), \beta =(01) , \mathrm{Tr}_\beta = M, q^c= 1$. 
\medskip

\noindent For $n=2$ we have
\medskip

\centerline{
\begin{tabular}{|c|c|c|c|c|c|c|c|c|}
\hline
$\pi$&$2-|\pi|$&$\pi\alpha$&$\beta$&$\mathrm{Tr}_\beta$&$l$&$q^{c}$\\
\hline
id &$2$&$(012)$&$(012)$&$M^2$&0&1\\
\hline
$(12)$&$1$&$(02)(1)$&(0)&$1$&$1$&$1$\\
\hline
\end{tabular}}

\medskip

\noindent For $n=3$ 
\medskip

\centerline{
\begin{tabular}{|c|c|c|c|c|c|c|}
\hline
$\pi$&$3-|\pi|$&$\pi\alpha$&$\beta$&$\mathrm{Tr}_\beta $&$l$&$q^{c}$\\
\hline
id &$3$&$(0123)$&$(0123)$&$M^3$&$0$&$1$\\
\hline
$(1)(23)$&$2$&$(013)(2)$&$(01)$&$M$&$1$&$1$\\
\hline
$(2)(13)$&$2$&$(03)(12)$&$(0)(1)$&$\Tr M$&$1$&$N^{-1}$\\
\hline
$(12)(3)$&$2$&$(023)(1)$&$(01)$&$M$&$1$&$1$\\
\hline
\end{tabular}
}
\medskip

\noindent For $n=4$
\medskip

\centerline{
\begin{tabular}{|c|c|c|c|c|c|c|}
\hline
$\pi$&$4 - |\pi|$&$\pi\alpha$&$\beta$&$\mathrm{Tr}_\beta$&$l$&
$q^c
$\\
\hline
id &$4$&$(01234)$&$(01234)$&$M^4$&$0$&$1$\\
\hline
$(1)(23)(4)$&$3$&$(0134)(2)$&$(012)$&$M^2$&$1$&$1$\\
\hline
$(1)(2)(34)$&$3$&$(0124)(3)$&$(012)$&$M^2$&$1$&$1$\\
\hline
$(1)(24)(3)$&$3$&$(014)(23)$&$(01)(2)$&$M\Tr M$&$1$&$N^{-1}$\\
\hline
$(14)(2)(3)$&$3$&$(04)(123)$&$(0)(12)$&$\Tr M^2$&$1$&$N^{-1}$\\
\hline
$(13)(2)(4)$&$3$&$(034)(12)$&$(02)(1)$&$M \Tr M$&$1$&$N^{-1}$\\
\hline
$(12)(3)(4)$&$3$&$(0234)(1)$&$(012)$&$M^2$&$1$&$1$\\
\hline
$(12)(34)$&$2$&$(024)(1)(3)$&$(0)$&$1$&$2$&$1$\\
\hline
$(13)(24)$&$2$&$(03214)$&$(0)$&$1$&$2$&$N^{-2}$\\
\hline
$(14)(23)$&$2$&$(04)(13)(2)$&$(0)$&$1$&$2$&$1$\\
\hline
\end{tabular}
}
\subsection{Proof of Proposition \ref{prop:ab3}}

We use formulas (\ref{defa}) and (\ref{defb}).

\noindent For $n=2$, we have  $\tilde \alpha= (0)(12)(34)$, $\pi \in \mathcal P_2(2,2)$
\medskip

\centerline{
\begin{tabular}{|c|c|c|c|c|}
\hline
$\pi$&$\pi\tilde\alpha$&$\cyc_0$&$\mathcal P'$\\ 
\hline
$(13)(24)$&$(0)(14)(23)$&$2$&$0$\\
\hline
$(14)(23)$&$(0)(13)(24)$&$2$&$1$\\
\hline
\end{tabular}
}
\smallskip

\noindent For $n=3$, we have  $\tilde\alpha =(0)(123)(456), \pi\in \mathcal P_2(3,3)$
\medskip

\centerline{
\begin{tabular}{|c|c|c|c|}
\hline
$\pi$&$\pi\tilde \alpha$&$\cyc_0$&$
\mathcal P'
$\\
\hline
$ (14)(25)(36)$&$(0) (153426)$&$1$&$1$\\
\hline
$(14)(26)(35)$&$(0) (16)(25)(34)$&$3$&$0$\\
\hline
$(15)(24)(36)$&$(0) (14)(26)(35)$&$3$&$0$\\
\hline
$ (15)(26)(34)$&$(0) (163524)$&$1$&$1$\\
\hline
$ (16)(24)(35)$&$(0) (143625)$&$1$&$1$\\
\hline
$(16)(25)(34)$&$(0) (15)(24)(36)$&$3$&$1$\\
\hline
\end{tabular}
} 
\subsection{Proof of Proposition \ref{prop:sigmahat}}
We apply formula (\ref{additional}).

\medskip

1) $\hat\sigma= (0)(12)$

\centerline{
\begin{tabular}{|c|c|c|c|c|c|c|}
\hline
$\pi$&$2-|\pi|$&$\pi\alpha$&$\beta$&$\mathrm{Tr}_\beta $&$l$&$q^{c}$\\
\hline
\hbox{id}&$2$&$(0)(12)$&$(0)(12)$&$\Tr (X^2)$&$0$&$1$\\
\hline
$(12)$&$1$&$(0)(1)(2)$&$(0)$&$1$&$1$&$N$\\
\hline
\end{tabular}
}
\medskip

2) $\hat\sigma= (0)(1)(2)$

\centerline{
\begin{tabular}{|c|c|c|c|c|c|c|}
\hline
$\pi$&$2-|\pi|$&$\pi\alpha$&$\beta$&$\mathrm{Tr}_\beta $&$l$&$q^{c}$\\
\hline
\hbox{id}&$2$&$\mathrm{id}$&$\mathrm{id}$&$(\Tr X)^2$&$0$&$1$\\
\hline
$(12)$&$1$&$(0)(12)$&$(0)$&$1$&$1$&$1$\\
\hline
\end{tabular}
}
\medskip

3) $\hat\sigma = (0)(123)$

\centerline{
\begin{tabular}{|c|c|c|c|c|c|c|}
\hline
$\pi$&$3-|\pi|$&$\pi\alpha$&$\beta$&$\Tr_\beta $&$l$&$q^{c}$\\
\hline
\hbox{id}&$3$&$(0)(123)$&$(0)(123)$&$\Tr (X^3)$&$0$&$1$\\
\hline
$(12)(3)$&$2$&$(0)(1)(23)$&$(0)(1)$&$\Tr X$&$1$&$1$\\
\hline
$(13)(2)$&$2$&$(0)(12)(3)$&$(0)(1)$&$\Tr X$&$1$&$1$ \\
\hline
$(1)(23)$&$2$&$(0)(13)(2)$&$(0)(1)$&$\Tr X$&$1$&$1$\\
\hline
\end{tabular}
}
\medskip

4) $\hat\sigma = (0)(12)(3)$

\centerline{
\begin{tabular}{|c|c|c|c|c|c|c|}
\hline
$\pi$&$3-|\pi|$&$\pi\alpha$&$\beta$&$\Tr_\beta $&$l$&$q^{c}$\\
\hline
\hbox{id}&$3$&$(0)(12)(3)$&$(0)(12)(3)$&$\Tr (X^2)\Tr X$&$0$&$1$\\
\hline
$(12)(3)$&$2$&$\hbox{id}$&$(0)(1)$&$\Tr X$&$1$&$N$\\
\hline
$(1)(23)$&$2$&$(0)(132)$&$(0)(1)$&$\Tr X$&$1$&$N^{-1}$\\
\hline
$(13) (2)$&2$$&$(0)(123)$&$(0)(1)$&$\Tr X$&$1$&$N^{-1}$\\
\hline
\end{tabular}
}
\medskip

5) $\hat\sigma =(0)(1)(2)(3)$

\centerline{
\begin{tabular}{|c|c|c|c|c|c|c|}
\hline
$\pi$&$3-|\pi|$&$\pi\alpha$&$\beta$&$\Tr_\beta $&$l$&$q^{c}$\\
\hline
\hbox{id}&$3$&$\hbox{id}$&$\hbox{id}$&$(\Tr X)^3$&$0$&$1$\\
\hline
$(12)(3)$&$2$&$(0)(12)(3)$&$(0)(1)$&$\Tr X$&$1$&$1$\\
\hline
$(13)(2)$&$2$&$(0)(13)(2)$&$(0)(1)$&$\Tr X$&$1$&$1$\\
\hline
$(1)(23))$&$2$&$(0)(1)(23)$&$(0)(1)$&$\Tr X$&$1$&$1$\\
\hline
\end{tabular}
}

\section*{Acknowledgement} 
A.R. would like to thank Jos\' e Le\'on for his valuable discussions during the preparation of the article \cite{leon2021} on large deviations and for providing  us with Notarnicola's reference.
\bibliographystyle{plain}
\bibliography{bibfree}
\end{document}